\documentclass[11pt, reqno]{amsart}

 \usepackage{amsmath}
 \usepackage{amssymb}
\usepackage{bm}

\newcommand{\mysection}[1]{\section{#1}
      \setcounter{equation}{0}}

\newcommand{\nlimsup}{\operatornamewithlimits{\overline{lim}}}
\newcommand{\nliminf}{\operatornamewithlimits{\underline{lim}}}

\newtheorem{theorem}{Theorem}[section]
\newtheorem{lemma}[theorem]{Lemma}

\theoremstyle{definition}
\newtheorem{assumption}{Assumption}[section]
\newtheorem{definition}{Definition}[section]
\theoremstyle{remark}
\newtheorem{remark}{Remark}[section]

\newcommand{\tr}{\text{\rm tr}\,}

\newcommand\bbeta{\text{\raise-.2ex\hbox{$\bm{\beta}$}}}

\newcommand\bR{\mathbb{R}}

\newcommand\bB{\mathbb{B}}

\newcommand\cF{\mathcal{F}}

\newcommand\cL{\mathcal{L}}

\newcommand\frA{\mathfrak{A}}
\newcommand\frB{\mathfrak{B}}

\newcommand\infsup{\operatornamewithlimits{inf\,\,\,sup}}
\newcommand\supinf{\operatornamewithlimits{sup\,\,\,inf}}
\newcommand\infinf{\operatornamewithlimits{inf\,\,\,inf}}

\begin{document}

\title[Independence of the value function]
{On the independence of the value function
for stochastic differential games of
the probability space}

\author{N.V. Krylov}
\thanks{The  author was partially supported by
 NSF Grant DMS-1160569}
\email{krylov@math.umn.edu}
\address{127 Vincent Hall, University of Minnesota,
 Minneapolis, MN, 55455}

\keywords{Stochastic differential games,
Isaacs equation, value functions}

\subjclass[2010]{49N70, 35D40, 49L25}

\begin{abstract}
We show that the value function in a stochastic differential
game
does not change  if we keep the same 
space $(\Omega,\cF)$ but introduce probability measures
by means of  Girsanov's transformation {\em depending\/}
on the policies of the players. We also
show that the value function
does not change if we allow the driving Wiener processes
to depend on the policies of the players. Finally,
we  
show that the value function
does not change if we perform a random time change
with the rate depending on the policies of the players.
 
\end{abstract}

\maketitle

\mysection{Introduction}

Let $\bR^{d}=\{x=(x^{1},...,x^{d})\}$
be a $d$-dimensional Euclidean space and  let $d_{1}\geq d$
  be an integer.
Assume that we are given separable metric spaces
  $A$ and $B$,   and let,
for each $\alpha\in A$, $\beta\in B$, 
  the following 
  functions on $\bR^{d}$ be given: 

(i) $d\times d_{1}$
matrix-valued $\sigma^{\alpha\beta}( x)
=\sigma(\alpha,\beta, x)=
(\sigma^{\alpha\beta}_{ij}( x))$,

(ii)
$\bR^{d}$-valued $b^{\alpha\beta}( x)=
b(\alpha,\beta, x)=
(b^{\alpha\beta}_{i }( x))$, and

(iii)
real-valued  functions $c^{\alpha\beta}( x)=c(\alpha,\beta, x)\geq0$,   
  $f^{\alpha\beta}( x)=f(\alpha,\beta, x)$, and  
$g(x)$. 
 
Under natural assumptions which will be specified later,
on a probability space
$(\Omega,\cF,P)$ carrying a $d_{1}$-dimensional Wiener process
$w_{t}$
one   associates with these objects and a bounded domain 
$G\subset\bR^{d}$
a stochastic differential
game 
with the diffusion term $\sigma^{\alpha\beta}(x)$,
  drift term $b^{\alpha\beta}(x)$, discount rate 
$c^{\alpha\beta}(x)$, running cost $f^{\alpha\beta}(x)$,
and the final cost $g(x)$ payed when the underlying process
first exits from $G$.

After the order of players is specified in a certain way 
it turns 
out (see our Remark \ref{remark 3.14.1}) that the value function $v(x)$
of this differential game is a unique 
continuous in $\bar{G}$ viscosity
solution
of the Isaacs equation
\begin{equation}
                                                      \label{3.14.3}
H[v]=0
\end{equation}
in $G$ with boundary condition $v=g$ on $\partial G$,
where for a sufficiently smooth function $u=u(x)$  
\begin{equation}
                                                     \label{1.16.1}
H[u](x)=\supinf_{\alpha\in A\,\,\beta\in B}
[L^{\alpha\beta} u(x)+f^{\alpha\beta} (x)],
\end{equation}
$$
L^{\alpha\beta} u(x):=a^{\alpha\beta}_{ij}( x)D_{ij}u(x)+
b ^{\alpha\beta}_{i }( x)D_{i}u(x)-c^{\alpha\beta} ( x)u(x),
$$
$$
a^{\alpha\beta}( x):=(1/2)\sigma^{\alpha\beta}( x)
(\sigma^{\alpha\beta}( x))^{*},\quad
D_{i}=\partial/\partial x^{i},\quad D_{ij}=D_{i}D_{j}.
 $$

We will assume that $\sigma$ and $b$ are uniformly Lipschitz
with respect to $x$, $\sigma\sigma^{*}$
is uniformly nondegenerate, and $c$ and $f$ are uniformly bounded.
In such a situation
uniqueness of continuous viscosity solutions or even
continuous $L_{p}$ viscosity solutions
 of \eqref{1.16.1} is shown in \cite{JS05}
and therefore the fact of the independence
of $v$ of the probability space seems to be obvious.

Roughly speaking, the goal of this paper is to show that the value function
does not change even if we keep the same 
space $(\Omega,\cF)$ but introduce probability measures
by means of  Girsanov's transformation {\em depending\/}
on the policies of the players. We also
show that the value function
does not change if we allow the driving Wiener processes
to depend on the policies of the players. Finally,
we  
show that the value function
does not change if we perform a random time change
with the rate depending on the policies of the players.

These facts are well known for controlled diffusion
processes and play there a very important role,
in particular, while estimating the derivatives
of the value function. A rather awkward substitute
of them for stochastic differential games
was used for the same purposes in \cite{Kr_14_1}.
Applying the results presented here one can
make many constructions in \cite{Kr_14_1}
more natural and avoid introducing
auxiliary ``shadow'' processes.

However, not all proofs
in \cite{Kr_14_1} can be simplified
using our present methods. We deliberately
avoided discussing the way to use
the external parameters in contrast with \cite{Kr_14_1}
just to make the presentation more transparent.

Our proofs do not use anything from the theory
of viscosity solutions and are  based on a version of
\'Swi{\c e}ch's (\cite{Sw96}) idea as presented 
in \cite{Kr_13_1} and a general solvability theorem
in class $C^{1,1}$ of Isaacs equations from \cite{Kr_12}.

The article is organized as follows. In Section
\ref{section 2.26.3} we present our main result, Theorem \ref{theorem 1.14.1}.
We prove it in Section \ref{section 2.11.1} under the additional
assumption
that the corresponding Isaacs equation has a smooth solution.
Then in Section \ref{section 4.6.1} we   allow the
solutions to belong to the Sobolev class $W^{2}_{d}$.
 Section \ref{section 2.25.1} contains a general approximation
result, which allows us in Section \ref{section 4.6.2} to
use a result from \cite{Kr_12} (see Theorem \ref{theorem 9.23.1})
and conclude the proof of Theorem \ref{theorem 1.14.1}
in the general case.
 
\mysection{Main result }
                                           \label{section 2.26.3}

We start with our assumptions.
 
\begin{assumption}
                                    \label{assumption 1.9.1}

(i) The functions $\sigma^{\alpha\beta}(x)$, $b^{\alpha\beta}(x)$,
$c^{\alpha\beta}(x)$, and $f^{\alpha\beta}(x)$
are continuous with respect to
$\beta\in B$ for each $(\alpha, x)$ and continuous with respect
to $\alpha\in A$ uniformly with respect to $\beta\in B$
for each $x$. The function $g(x)$ is bounded and continuous.  

(ii) The functions 
$c^{\alpha\beta}( x )$
and
$f^{\alpha\beta}( x )$
are uniformly continuous with
respect to $x$ uniformly with respect to 
$(\alpha,\beta )\in A\times B $
and
for any $x \in\bR^{d}$ 
 and $(\alpha,\beta )\in A\times B $
$$
  \|\sigma^{\alpha\beta}( x )\|,|b^{\alpha\beta}( x )|,
|c^{\alpha\beta}( x )|,|f^{\alpha\beta}( x )|
\leq K_{0},
$$
where   $K_{0}$ is
a fixed constants and for a matrix $\sigma$ we denote $\|\sigma\|^{2}
=\tr\sigma\sigma^{*}$,

(iii) For any $(\alpha,\beta )\in A\times B $ and $x,y\in\bR^{d}$
we have
$$
\|\sigma^{\alpha\beta}( x )-\sigma^{\alpha\beta}( y )\|+
|b^{\alpha\beta}( x )-b^{\alpha\beta}( y )|\leq K_{0}|x-y|.
$$
\end{assumption}

Let $(\Omega,\cF,P)$ be a complete probability space,
let $\{\cF_{t},t\geq0\}$ be an increasing filtration  
of $\sigma$-fields $\cF_{t}\subset \cF $ such that
each $\cF_{t}$ is complete with respect to $\cF,P$.

 The set of progressively measurable $A$-valued
processes $\alpha_{t}=\alpha_{t}(\omega)$ is denoted by $\frA$. 
Similarly we define $\frB$
as the set of $B$-valued  progressively measurable functions.
By  $ \bB $ we denote
the set of $\frB$-valued functions 
$ \bbeta(\alpha_{\cdot})$ on $\frA$
such that, for any $T\in(0,\infty)$ and any $\alpha^{1}_{\cdot},
\alpha^{2}_{\cdot}\in\frA$ satisfying
\begin{equation}
                                                  \label{4.5.4}
P(  \alpha^{1}_{t}=\alpha^{2}_{t} 
 \quad\text{for almost all}\quad t\leq T)=1,
\end{equation}
we have
$$
P(  \bbeta_{t}(\alpha^{1}_{\cdot})=\bbeta_{t}(\alpha^{2}_{\cdot}) 
\quad\text{for almost all}\quad t\leq T)=1.
$$
 
\begin{definition}
A 
function $p^{\alpha_{\cdot}\beta_{\cdot}}_{t}=
p^{\alpha_{\cdot}\beta_{\cdot}}_{t}(\omega)$ given
on $\frA\times\frB\times\Omega\times[0,\infty)$
with values in some measurable space
is called {\em a control adapted process} if,
for any $(\alpha_{\cdot},\beta_{\cdot})\in\frA\times\frB$,
it is progressively measurable in $(\omega,t)$ and, for any
$T\in(0,\infty)$, we have
$$
P( p^{\alpha^{1}_{\cdot}\beta^{1}_{\cdot}}_{t}
=p^{\alpha^{2}_{\cdot}\beta^{2}_{\cdot}}_{t}
\quad\text{for almost all}\quad t\leq T)=1
$$
as long as
$$
P(  \alpha^{1}_{t}=\alpha^{2}_{t} ,
\beta^{1}_{t}=\beta^{2}_{t}
\quad\text{for almost all}\quad t\leq T)=1.
$$
 
\end{definition}

\begin{assumption}
                                                \label{assumption 2.9.1}
For each $\alpha_{\cdot}\in\frA$ and $\beta_{\cdot}\in\frB$
we are given control adapted processes

(i) $w^{\alpha_{\cdot}\beta_{\cdot}}_{t},t\geq0$, which
are  standard $d_{1}$-dimensional Wiener process relative to
  to the filtration $\{\cF_{t},t\geq0\}$,

(ii)  $r^{\alpha_{\cdot}\beta_{\cdot}}_{t},t\geq0$,
and $\pi^{\alpha_{\cdot}\beta_{\cdot}}_{t},t\geq0$, which
are real-valued and $\bR^{d_{1}}$-valued, respectively,

(iii) for all values of the arguments
$$
\delta_{1}^{-1}\geq r^{\alpha_{\cdot}\beta_{\cdot}}_{t}\geq\delta_{1},\quad
|\pi^{\alpha_{\cdot}\beta_{\cdot}}_{t}|\leq K_{1},
$$
where $\delta_{1}>0$ and $K_{1}\in(0,\infty)$ are fixed constants.
\end{assumption}

Finally we 
introduce
$$
a^{\alpha\beta}( x):=(1/2)\sigma^{\alpha\beta}( x)
(\sigma^{\alpha\beta}( x))^{*},
$$
fix a domain $G\subset\bR^{d}$, and impose the following.

\begin{assumption}
                                           \label{assumption 3.19.1}
 $G$ is a bounded domain of class $C^{2}$
and there exists a constant $\delta\in(0,1)$ such that
for any $\alpha\in A$, $\beta\in B$,  
and $x,\lambda\in \bR^{d}$
$$
\delta|\lambda|^{2}\leq a^{\alpha\beta}_{ij}( x)
\lambda^{i}\lambda^{j}\leq \delta^{-1}|\lambda|^{2}.
$$

\end{assumption}

\begin{remark}
                                                     \label{remark 2.22.2}

As is well known, if Assumption \ref{assumption 3.19.1} 
is satisfied, then
 there exists a bounded from above
  $\Psi\in   C^{2}_{loc}(\bR^{d})$
such that 
  $\Psi>0$ in $G$, $\Psi=  0$ on $\partial G$, and
for all $\alpha\in A$, $\beta\in B$,  
and $x\in G$
\begin{equation}
                                       \label{3.20.1}
  L^{\alpha\beta}\Psi( x)
+c^{\alpha\beta}\Psi( x)\leq-1.
\end{equation}
\end{remark}

For $\alpha_{\cdot}\in\frA$, 
$ \beta_{\cdot} 
\in\frB$, and $x\in\bR^{d}$ consider the following It\^o
equation
$$
x_{t}=x+\int_{0}^{t}
r^{\alpha_{\cdot}\beta_{\cdot}}_{s}
\sigma^{\alpha_{s}\beta_{s} }(  x_{s})\,dw^{\alpha_{\cdot}\beta_{\cdot}}_{s}
$$
\begin{equation}
                                             \label{5.11.1}
+\int_{0}^{t}[r^{\alpha_{\cdot}\beta_{\cdot}}_{s}]^{2}\big[b^{\alpha_{s}
\beta_{s} }(  x_{s})+\sigma^{\alpha_{s}\beta_{s} }(  x_{s})
\pi^{\alpha_{\cdot}\beta_{\cdot}}_{s}\big]\,ds.
\end{equation}

Observe that equation \eqref{5.11.1} satisfies the usual hypothesis, that is
for any $\alpha_{\cdot}\in\frA$, 
$ \beta_{\cdot} 
\in\frB$,  $x\in\bR^{d}$, and $T\in(0,\infty)$ it has a unique solution
on $[0,T]$
denoted by $x^{\alpha_{\cdot} 
\beta_{\cdot} x}_{t} $ and $x^{\alpha_{\cdot} 
\beta_{\cdot} x}_{t} $ is a control adapted process
for each $x$.

Set
$$
\phi^{\alpha_{\cdot}\beta_{\cdot} x}_{t}
=\int_{0}^{t}[r^{\alpha_{\cdot}\beta_{\cdot}}_{s}]^{2}
c^{\alpha_{s}
\beta_{s} }( x^{\alpha_{\cdot}\beta_{\cdot}  x}_{s}) 
\,ds,
$$
$$
\psi^{\alpha_{\cdot}\beta_{\cdot} x}_{t}
=-(1/2)
\int_{0}^{t}[r^{\alpha_{\cdot}\beta_{\cdot}}_{s}]^{2}
|\pi^{\alpha_{\cdot}\beta_{\cdot}}_{s}|^{2}
\,ds-\int_{0}^{t}
r^{\alpha_{\cdot}\beta_{\cdot}}_{s}
\pi^{\alpha_{\cdot}\beta_{\cdot}}_{s}\,dw^{\alpha_{\cdot}\beta_{\cdot}}_{s},
$$
 define $\tau^{\alpha_{\cdot}\beta_{\cdot} x}$ as the first exit
time of $x^{\alpha_{\cdot} 
\beta_{\cdot} x}_{t}$ from $G$, and introduce
\begin{equation}
                                                    \label{2.12.2}
v(x)=\infsup_{\bbeta\in\bB\,\,\alpha_{\cdot}\in\frA}
E_{x}^{\alpha_{\cdot}\bbeta(\alpha_{\cdot})}\big[\int_{0}^{\tau}
r_{t}^{2}f( x_{t})e^{-\phi_{t}
-\psi_{t}}\,dt+g(x_{\tau})e^{-\phi_{\tau}
-\psi_{\tau}}\big],
\end{equation}
where the indices $\alpha_{\cdot}$, $\bbeta$, and $x$
at the expectation sign are written  to mean that
they should be placed inside the expectation sign
wherever and as appropriate, that is
$$
E_{x}^{\alpha_{\cdot}\beta_{\cdot}}\big[\int_{0}^{\tau}r^{2}_{t}
f( x_{t})e^{-\phi_{t}-\psi_{t}}\,dt+g(x_{\tau})e^{-\phi_{\tau}
-\psi_{\tau}}\big]
$$
$$
:=
E \big[
 g(x^{\alpha_{\cdot}\beta_{\cdot}  x}
_{\tau^{\alpha_{\cdot}\beta_{\cdot}  x}}
)
e^{-\phi^{\alpha_{\cdot}\beta_{\cdot}  x}
_{\tau^{\alpha_{\cdot}\beta_{\cdot}  x}}
-\psi^{\alpha_{\cdot}\beta_{\cdot}  x}
_{\tau^{\alpha_{\cdot}\beta_{\cdot}  x}}}
$$
$$
+\int_{0}^{\tau^{\alpha_{\cdot}\beta_{\cdot}  x}}
[r^{\alpha_{\cdot}\beta_{\cdot}  x}_{t}]^{2}
f^{\alpha_{t}\beta_{t}   }
( 
x^{\alpha_{\cdot}\beta_{\cdot}  x}_{t})
e^{-\phi^{\alpha_{\cdot}\beta_{\cdot}  x}_{t}
-\psi^{\alpha_{\cdot}\beta_{\cdot}  x}_{t}}\,dt\big].
$$
Observe that, formally, the value $x_{\tau}$ 
 may not be defined if $\tau=\infty$.
In that case we set the corresponding terms to equal zero.
 The above definitions make  perfect sense due to our Remark \ref{remark 2.20.1}.

Here is our   main result. 

\begin{theorem}
                                                    \label{theorem 1.14.1}
 
Under the above assumptions 
the function $v(x)$  
is independent of the choice of the probability space,
filtration and control adapted process 
$(r,\pi,w)^{\alpha_{\cdot}\beta_{\cdot}}_{t}$,
it is bounded and continuous in $\bar{G}$.
\end{theorem}

\begin{remark}
                                                       \label{remark 3.14.1}
Once we know that $v(x)$  
is independent of the choice of the probability space,
filtration and control adapted process 
$(r,\pi,w)^{\alpha_{\cdot}\beta_{\cdot}}_{t}$, we can take any
probability space carrying a $d_{1}$-dimensional Wiener process
$w_{t}$ and construct $v(x)$ by setting $w^{\alpha_{\cdot}\beta_{\cdot}}_{t}
=w_{t}$, $r\equiv1$, $\pi\equiv0$. In that case we are
in the position to apply the results of 
  \cite{FS89}, \cite{Ko09},  and \cite{Kr_13}
according to which $v$ is continuous in $\bar{G}$
and satisfies the dynamic programming principle.
Then it is a standard fact that $v$ is a viscosity
solution of \eqref{3.14.3} (see, for instance, \cite{FS89}, \cite{Ko09},
\cite{Sw96}). Indeed, if a smooth function $\psi(x)$ is such that
$\psi(x)\geq v(x)$ in a neighborhood of $x_{0}\in G$ and
$\psi(x_{0})= v(x_{0})$, then by defining 
$\gamma^{\alpha_{\cdot}\beta_{\cdot}}_{\varepsilon}$,
$\varepsilon>0$ as the first exit time of $x^{\alpha_{\cdot}\beta_{\cdot}x_{0}}_{t}$
from an $\varepsilon$-neighborhood of $x_{0}$ for all small
$\varepsilon$ we have
$$
\psi(x_{0})=v(x_{0})
=\infsup_{\bbeta\in\bB\,\,\alpha_{\cdot}\in\frA}
E_{x_{0}}^{\alpha_{\cdot}\bbeta(\alpha_{\cdot})}
\big[\int_{0}^{\gamma_{\varepsilon}}
 f( x_{t})e^{-\phi_{t}}\,dt+v(x_{\gamma_{\varepsilon}})
e^{-\phi_{\gamma_{\varepsilon}}}\big]
$$
$$
\leq \infsup_{\bbeta\in\bB\,\,\alpha_{\cdot}\in\frA}
E_{x_{0}}^{\alpha_{\cdot}\bbeta(\alpha_{\cdot})}
\big[\int_{0}^{\gamma_{\varepsilon}}
 f( x_{t})e^{-\phi_{t}}\,dt+\psi(x_{\gamma_{\varepsilon}})
e^{-\phi_{\gamma_{\varepsilon}}}\big].
$$
On the other hand set $H[\psi]=-h$ and observe 
that by Theorem 4.1 of \cite{Kr_13_1}
$$
\psi(x_{0})=
\infsup_{\bbeta\in\bB\,\,\alpha_{\cdot}\in\frA}
E_{x_{0}}^{\alpha_{\cdot}\bbeta(\alpha_{\cdot})}
\big[\int_{0}^{\gamma_{\varepsilon}}
 (f+h)( x_{t})e^{-\phi_{t}}\,dt+\psi(x_{\gamma_{\varepsilon}})
e^{-\phi_{\gamma_{\varepsilon}}}\big].
$$
It follows that
\begin{equation}
                                                     \label{3.14.4}
\infinf_{\alpha_{\cdot}\in\frA\,\,\,\beta_{\cdot}\in\frB}
E_{x_{0}}^{\alpha_{\cdot}\beta _{\cdot}}
 \int_{0}^{\gamma_{\varepsilon}}
 h( x_{t})e^{-\phi_{t}}\,dt\leq0, 
\end{equation}
and if we assume that $H[\psi](x_{0})<0$, then $h>0$
in an $\varepsilon$-neighborhood of $x_{0}$
and \eqref{3.14.4} is impossible, since $c$ is bounded
and $\sigma$ and $b$ are bounded so that 
$E_{x_{0}}^{\alpha_{\cdot}\beta _{\cdot}}
\gamma_{\varepsilon}$ is bounded away from zero.
Hence $H[\psi](x_{0})\geq0$ and $v$ is a viscosity subsolution
by definition. Similarly one shows that
it is a viscosity supersolution.

Provided that we know that continuous viscosity solutions 
are unique the above argument proves the fact that the value function is
independent of the probability space (if we drop out $r$ and $\pi$
and take $w$ independent of the policies).
Jensen \cite{Je_88} proved uniqueness 
for {\em Lipschitz\/} continuous viscosity
solutions to the fully nonlinear second 
order elliptic PDE not explicitly depending on
$x$ on a bounded domain. Related results in the same year 
with $H$ depending on $x$
were published in Jensen-Lions-Souganidis
 \cite{JLS_88}.

In what concerns uniformly
nondegenerate Isaacs
equations, Trudinger in \cite{Tr_89} proves the existence and uniqueness
 of continuous viscosity solutions for Isaacs equations
if the coefficients are continuous and $a$
is $1/2$ H\"older continuous uniformly with respect to $\alpha,\beta$
(see Corollary 3.4 there).
Uniqueness   is also stated for Isaacs equations
with Lipschitz continuous $a$
as Corollary 5.11 in \cite{Je_98}.
 Jensen and \'Swi{\c e}ch
in \cite{JS05} further relaxed the requirement on $a$
and proved uniqueness of  continuous even $L_{p}$-viscosity
solutions.

\end{remark}

We will use Theorem \ref{theorem 1.14.1} to prove in a subsequent article
a result to state which
we need a few new objects.  
In the end of Section 1 of \cite{Kr_12} a 
function $P(u_{ij},u_{i},u)$ is constructed
defined for all symmetric $d\times d$ matrices $(u_{ij})$,
$\bR^{d}$-vectors $(u_{i})$, and $u\in\bR$ such that
  it is positive-homogeneous of degree
one, is
Lipschitz continuous, and at all points of differentiability of $P$
for all values of arguments we have $P_{u}\leq0$ and
$$
\hat\delta|\lambda|^{2}\leq P_{u_{ij}}
\lambda^{i}\lambda^{j}\leq \hat{\delta}^{-1}|\lambda|^{2},
$$
where $\hat{\delta}$ is a constant in $(0,1)$
depending only on
$d,K_{0}$, and $\delta$. For smooth enough functions $u(x)$ introduce
$$
P[u](x)=P(D_{ij}u(x),D_{i}u(x),u(x))
$$

We now state part of Theorem 1.1 of \cite{Kr_12}
which we need even in the present article.
\begin{theorem}
                                    \label{theorem 9.23.1}
Let $g\in C^{1,1}(\bR^{d})$. Then
for any $K\geq0$   the equation
\begin{equation}
                                               \label{9.23.2}
\max(H[u],P[u]-K)=0
\end{equation}
in $G$  (a.e.) with boundary condition $u=g$ on $\partial G$
has a unique solution 
$u\in C^{0,1}(\bar{G})\cap C^{1,1}_{loc}(G)$.
\end{theorem}

The result we are aiming at in a subsequent
article consists of proving
the conjecture stated in \cite{Kr_12}:

\begin{theorem}
                                                    \label{theorem 1.14.01}
 
Denote by $u_{K}$ the function from Theorem \ref{theorem 9.23.1}
and assume that $G$ and $g$ are of class $C^{3}$.
Then there exists a constant $N$
such that
$|v-u_{K}|\leq N/K$ on $G$ for $K\geq1$.
\end{theorem}
A very week version of this theorem was already used
in \cite{Kr_14_2} for establishing
a rate of convergence of finite-difference
approximations for solutions of Isaacs equations.

We finish this section with a useful technical result.

\begin{lemma}
                                 \label{lemma 2.20.1}
 
For any $\alpha_{\cdot}\in\frA$, $\beta_{\cdot}
\in\frB$, and $x\in\bR^{d}$ the process
$$
\exp(-\psi_{t\wedge\tau^{\alpha_{\cdot}\beta_{\cdot}x}}
^{\alpha_{\cdot}\beta_{\cdot}x})
$$
is a uniformly integrable martingale
on $[0,\infty)$. Furthermore, there exists
a constant $N$ independent of 
$\alpha_{\cdot}\in\frA$, $\beta_{\cdot}
\in\frB$, and $x\in\bR^{d}$ such that
\begin{equation}
                                       \label{2.20.1}
E^{\alpha_{\cdot}\beta_{\cdot}}_{x}
\int_{0}^{
\tau}
e^{-\psi_{s}}\,ds\leq N.
\end{equation}
Finally,
\begin{equation}
                                       \label{2.20.2}
E^{\alpha_{\cdot}\beta_{\cdot}}_{x}
e^{-\psi_{\tau}}=1.
\end{equation}
\end{lemma}

Proof. Owing to \eqref{3.20.1} by It\^o's formula
$$
E^{\alpha_{\cdot}\beta_{\cdot}}_{x}
G(x_{t\wedge\tau})e^{-\psi_{t\wedge\tau}}
=G(x)+E^{\alpha_{\cdot}\beta_{\cdot}}_{x}
\int_{0}^{t\wedge\tau}
r^{2}_{s}[LG+cG](x_{s})e^{-\psi_{s}}\,ds
$$
$$
\leq
 G(x)-\delta_{1}^{2}
E^{\alpha_{\cdot}\beta_{\cdot}}_{x}\int_{0}^{
t\wedge\tau}
e^{-\psi_{s}}\,ds,
$$
and \eqref{2.20.1} follows. To prove \eqref{2.20.2}
use that
$$
1=E^{\alpha_{\cdot}\beta_{\cdot}}_{x}
e^{-\psi_{t\wedge\tau}}
=E^{\alpha_{\cdot}\beta_{\cdot}}_{x}
e^{-\psi_{\tau}}I_{\tau\leq t}
+E^{\alpha_{\cdot}\beta_{\cdot}}_{x}
e^{-\psi_{t}}I_{\tau> t},
$$
where the last term decreases as $t$ increases,
which is seen from the formula, and tends to zero as $t\to\infty$
since its integral with respect to $t$ over $[0,
\infty)$ is finite being equal to the 
left-hand side of \eqref{2.20.1}.

Finally, the first assertion of the lemma follows
from \eqref{2.20.2} due to the well-known properties
of martingales. The lemma is proved.

\begin{remark}
                                 \label{remark 2.20.1}
In light of the proof of Lemma \ref{lemma 2.20.1}
$$
0=\lim_{t\to\infty}E^{\alpha_{\cdot}\beta_{\cdot}}_{x}
e^{-\psi_{t}}I_{\tau> t}
=\lim_{t\to\infty}E^{\alpha_{\cdot}\beta_{\cdot}}_{x}
e^{-\psi_{\tau}}I_{\tau> t}=
E^{\alpha_{\cdot}\beta_{\cdot}}_{x}
e^{-\psi_{\tau}}I_{\tau=\infty}.
$$
Hence defining the terms containing $x_{\tau}$ as zero
on the set where $\tau=\infty$ is indeed natural.
Lemma \ref{lemma 2.20.1} shows that the function $v$ is well defined and
one can rewrite its definition   as
$$
v(x)=\infsup_{\bbeta\in\bB\,\,\alpha_{\cdot}\in\frA}
E_{x}^{\alpha_{\cdot}\bbeta(\alpha_{\cdot})}\big[\int_{0}^{\tau}
r_{t}^{2}f( x_{t})e^{-\phi_{t}
}\,dt+g(x_{\tau})\big]e^{-\psi_{\tau}},
$$
which calls for   changes of probability measure by using 
Girsanov's theorem.

\end{remark}

\mysection{Proof of Theorem \protect\ref{theorem 1.14.1}
in case that the Isaacs equation has a smooth solution}

                                          \label{section 2.11.1}

In this section we   replace Assumption \ref{assumption 1.9.1} (iii)
with a weaker one.

\begin{assumption}
                                             \label{assumption 2.16.1}
 
The functions 
$\sigma^{\alpha\beta}( x )$
and
$b^{\alpha\beta}( x )$
are uniformly continuous with
respect to $x$ uniformly with respect to 
$(\alpha,\beta )\in A\times B $.
\end{assumption}
However, this time there is no guarantee that
equation \eqref{5.11.1} has a unique solution
and we impose the following.

\begin{assumption}
                                           \label{assumption 5.23.2} 
Equation \eqref{5.11.1} satisfies the usual hypothesis, that is
for any $\alpha_{\cdot}\in\frA$, 
$ \beta_{\cdot} 
\in\frB$,  $x\in\bR^{d}$, and $T\in(0,\infty)$ it has a unique solution
on $[0,T]$
denoted by $x^{\alpha_{\cdot} 
\beta_{\cdot} x}_{t} $ and $x^{\alpha_{\cdot} 
\beta_{\cdot} x}_{t} $ is a control adapted process
for each $x$.

\end{assumption}

We also assume that we are given two  functions
 $\hat{u},\check{u} \in C^{2}(\bar{G})$.

\begin{theorem}
                                         \label{theorem 2.11.1}

(i) If $H[\hat{u}]\leq0$  
  in $G$ 
and $\hat{u}\geq g$ on $\partial G$, then $v\leq\hat{u}$ in $\bar{G}$.

(ii) If $H[\check{u}]\geq0$  
  in $G$ 
and $\check{u}\leq g$  on $\partial G$,
 then $v\geq\check{u}$ in $\bar{G}$.

(iii) If $\hat{u}$ and $\check{u}$ are as in (i) and (ii) and
  $\hat{u}=\check{u}$, then
$v$ is independent of the choice of the probability space,
filtration, $r,\pi$, and $w$.

 \end{theorem}

 We need three lemmas.
\begin{lemma}
                                                 \label{lemma 2.10.1}
Set   $\kappa_{n}(t)=[nt]/n$. Then there exists
a constant $N$ such that for all $n\geq1,x\in\bR^{d},
\alpha_{\cdot}\in\frA,\beta_{\cdot}\in\frB$
we have
\begin{equation}
                                                       \label{2.10.3}
E^{\alpha_{\cdot}\beta_{\cdot}}_{x}\int_{0}^{\tau}
e^{-\phi_{t}-\psi_{t}}|x_{t}-x_{\kappa_{n}(t)}|^{2}\,dt\leq N/n.
\end{equation}

\end{lemma}

Proof. For each fixed $t$ while estimating
$$
E^{\alpha_{\cdot}\beta_{\cdot}}_{x} 
e^{-\phi_{t}-\psi_{t}}
|x_{t}-x_{\kappa_{n}(t)}|^{2}I_{\tau>t}
$$
Girsanov's theorem allows us to assume that $\pi\equiv 0$.
In that case
for simplicity of notation we will drop the
indices $x,\alpha_{\cdot},\beta_{\cdot}$ and observe that
$$
|x_{t}-x_{\kappa_{n}(t)}|^{2}\leq 2\big|\int_{k/n}^{t}
\sigma(x_{s})\,dw_{s}\big|^{2}+
\big|\int_{k/n}^{t}
b(x_{s})\,ds\big|^{2}
$$
so that for $t\in[k/n,(k+1)/n]$
$$
E^{\alpha_{\cdot}\beta_{\cdot}}_{x}\{
e^{-\phi_{t} }|x_{t}-x_{k/n}|^{2}\mid\cF_{k/n}\}
\leq  E^{\alpha_{\cdot}\beta_{\cdot}}_{x}\{
 |x_{t}-x_{k/n}|^{2}\mid\cF_{k/n}\}\leq N/n,
$$
where $N$ depends only on $d$ and $K_{0}$. 
Hence, owing also to \eqref{2.20.1}
the left-hand side of  \eqref{2.10.3} is dominated by
$$
\int_{0}^{\infty}E^{\alpha_{\cdot}\beta_{\cdot}}_{x} 
e^{-\phi_{t}}|x_{t}-x_{\kappa_{n}(t)}|^{2}I_{\tau>t}\,dt
\leq \int_{0}^{\infty}E^{\alpha_{\cdot}\beta_{\cdot}}_{x} 
 |x_{t}-x_{\kappa_{n}(t)}|^{2}I_{\tau>
\kappa_{n}(t)}\,dt
$$
$$
\leq Nn^{-1} \int_{0}^{\infty}E^{\alpha_{\cdot}\beta_{\cdot}}_{x} 
 I_{\tau>
\kappa_{n}(t)}\,dt=
Nn^{-1}  E^{\alpha_{\cdot}\beta_{\cdot}}_{x} 
(\tau+1/n)\leq Nn^{-1}.
$$
The lemma is proved.

For a stopping time $\gamma$
we say that a process $\xi_{t}$ is a submartingale
on $[0,\gamma]$ if $\xi_{t\wedge\gamma}$ is a submartingale.
Similar definition applies to supermartingales.

The proof of the following lemma and Lemma \ref{lemma 1.17.2}
follows a
version of
\'Swi{\c e}ch's (\cite{Sw96}) idea as it is presented 
in \cite{Kr_13_1}.

\begin{lemma}
                                                  \label{lemma 1.16.1}
Let $H[\hat{u}]\leq0$  
  in $G$. Then for any $x\in\bR^{d}$,
 $\alpha_{\cdot}\in\frA$, and $\varepsilon>0$,
there exist a sequence $\beta^{n}_{\cdot}(\alpha_{\cdot})=
\beta^{n}_{\cdot}(\alpha_{\cdot},x,\varepsilon)
\in\frB$, $n=1,2,...$, and a sequence of increasing continuous
$\{\cF_{t}\}$-adapted processes 
$\eta^{n\varepsilon}_{t}(\alpha_{\cdot})
=\eta^{n\varepsilon}_{t}(\alpha_{\cdot},x)$
 with $\eta^{n\varepsilon}_{0}(\alpha_{\cdot})=0$
 such that 
\begin{equation}
                                                 \label{2.29.1}
\sup_{n}E 
\eta^{n\varepsilon}_{\infty}(\alpha_{\cdot})<\infty,
\end{equation}
 the processes
$$
\kappa^{n\varepsilon}_{t  }(\alpha_{\cdot}):=
\hat{u}(x^{n}_{t})
e^{-\phi^{n}_{t}-\psi^{n}_{t}}
-\eta^{n\varepsilon}_{t} (\alpha_{\cdot})
+\int_{0}^{t}[r^{n}_{t}]^{2}
f^{n}_{s} ( x^{n}_{s}) 
e^{-\phi^{n}_{s}-\psi^{n}_{s}}\,ds,
$$
where
\begin{equation}
                                                   \label{5.6.1}
(x^{n}_{t},\phi^{n}_{t},\psi^{n}_{t})
=(x_{t},\phi_{t},\psi_{t})
^{\alpha_{\cdot}\beta^{n}_{\cdot}(\alpha_{\cdot})  x},\quad
f^{n}_{t}( x)=f^{\alpha_{t}\beta^{n}_{t}(\alpha_{\cdot})}( x),\quad
r^{n}_{t}=r^{\alpha_{\cdot}\beta^{n}_{\cdot}(\alpha_{\cdot})}_{t}
\end{equation}
are supermartingales on $[0,\tau^{\alpha_{\cdot}\beta^{n}_{\cdot}(\alpha_{\cdot})
  x}]$, and
\begin{equation}
                                                       \label{3.3.1}
\nlimsup_{n\to\infty}
\sup_{\alpha_{\cdot}\in\frA}E 
 \eta^{n\varepsilon}_{\tau}(\alpha_{\cdot})\leq
 N \varepsilon  ,
\end{equation}
where $N$ is independent of $x$ and $\varepsilon$.
 Finally,
\begin{equation}
                                           \label{2.29.3}
\sup_{\alpha_{\cdot}\in\frA}\sup_{n}
E 
\sup_{t\geq0}|
\kappa^{n\varepsilon}_{t\wedge\tau }(\alpha_{\cdot})|<\infty .
\end{equation}

\end{lemma}  

Proof. Since $B$ is separable and $a^{\alpha\beta},b^{\alpha\beta},
c^{\alpha\beta}$, and $f^{\alpha\beta}$ are continuous with respect to
$\beta$ one can replace $B$ in \eqref{1.16.1} with an appropriate
countable subset $B_{0}=\{\beta_{1},\beta_{2},...\}$.
 Then for each $\alpha\in A$ and $x\in G$
define $\beta(\alpha,x)$ as $\beta_{i}\in B_{0}$ with the least 
$i$ such that
\begin{equation}
                                                  \label{1.17.2}
0
\geq L^{\alpha\beta_{i}} \hat{u}(x)+f^{\alpha\beta_{i}} 
(x)-\varepsilon.  
\end{equation}
For each $i$ the right-hand side of \eqref{1.17.2} is
Borel in $x$ and continuous in $\alpha$. Therefore,
it is a Borel function of $(\alpha,x)$, implying that  $\beta(\alpha,x)$ 
also is a Borel function
of $(\alpha,x)$. For $x\not\in G$ set $\beta(\alpha,x)=\beta^{*}$,
where $\beta^{*}$ is a fixed element of $B$.
Then we have that in $G$
\begin{equation}
                                                  \label{1.17.3}
0
\geq L^{\alpha\beta(\alpha,x)} \hat{u}(x)
+f^{\alpha\beta(\alpha,x)} (x)-\varepsilon.
\end{equation}

 After that fix $x$,
define  $\beta^{n0}_{t}(\alpha_{\cdot})=\beta(\alpha_{t},x)$, $t\geq0$,
and for $k\geq1$ introduce
$\beta^{nk}_{t}(\alpha_{\cdot})$ 
recursively so that
\begin{equation}
                                                  \label{4.5.2}
\beta^{nk}_{t}(\alpha_{\cdot})=\beta^{n(k-1)}_{t}(\alpha_{\cdot})
\quad\text{for}\quad t<k/n,
\end{equation}
$$
\beta^{nk}_{t}(\alpha_{\cdot})=\beta(\alpha_{t},x^{nk }_{k/n})
\quad\text{for}\quad t\geq k/n,
$$
where $x^{nk}_{t}$, $k=1,2,...$, is a unique solution of
$$
x _{t}=x+\int_{0}^{t}
r^{\alpha_{\cdot}\beta^{n(k-1)}_{\cdot}(\alpha_{\cdot})}_{s}
\sigma(\alpha_{s},\beta^{n(k-1)}_{s}(\alpha_{\cdot}),
 x_{s})\,
dw^{\alpha_{\cdot}\beta^{n(k-1)}_{\cdot}(\alpha_{\cdot})}_{s}
$$
\begin{equation}
                                                      \label{4.5.1}
+\int_{0}^{t}
[r^{\alpha_{\cdot}\beta^{n(k-1)}_{\cdot}} _{s}]^{2}
\big[b(\alpha_{s},\beta^{n(k-1)}_{s}(\alpha_{\cdot}),x_{s}
)+\sigma(\alpha_{s},\beta^{n(k-1)}_{s}(\alpha_{\cdot}),
 x_{s})\pi^{\alpha_{\cdot}\beta^{n(k-1)}_{\cdot}}_{s}
\big]\,
ds.
\end{equation}

To show that the above definitions make sense, observe that, 
by Assumption \ref{assumption 5.23.2},
$x^{n1}_{t}$ is well defined for all $t$.
 Therefore, $\beta^{n1}_{t}(\alpha_{\cdot})$
is also well defined, and by induction we conclude that
$x^{nk}_{t}$ and $\beta^{nk}_{t}(\alpha_{\cdot}) $ are
well defined for all $k\geq 1$.

Furthermore, owing to \eqref{4.5.2} it makes sense to define
$$
\beta^{n}_{t}(\alpha_{\cdot})=\beta^{nk}_{t}(\alpha_{\cdot})
\quad\text{for}\quad t<k/n.
$$
Notice that by definition 
 $x^{n}_{t}:=x^{\alpha_{\cdot}\beta^{n}_{\cdot}(\alpha_{\cdot})x}_{t}$
satisfies the equation
$$
x _{t}=x+\int_{0}^{t}
r^{\alpha_{\cdot}\beta^{n }_{\cdot}(\alpha_{\cdot})}_{s}
\sigma(\alpha_{s},\beta^{n }_{s}(\alpha_{\cdot}),
 x_{s})\,
dw^{\alpha_{\cdot}\beta^{n }_{\cdot}(\alpha_{\cdot})}_{s}
$$
\begin{equation}
                                                      \label{4.6.1}
+\int_{0}^{t}
[r^{\alpha_{\cdot}\beta^{n }_{\cdot}} _{s}]^{2}
\big[b(\alpha_{s},\beta^{n }_{s}(\alpha_{\cdot}),x_{s}
)+\sigma(\alpha_{s},\beta^{n }_{s}(\alpha_{\cdot}),
 x_{s})\pi^{\alpha_{\cdot}\beta^{n }_{\cdot}(\alpha_{\cdot})}_{s}
\big]\,
ds.
\end{equation}
For $t<k/n$ we have $\beta^{n}_{t}(\alpha_{\cdot})
=\beta^{n(k-1)}_{t}(\alpha_{\cdot})$,
so that for $t\leq k/n$ equation \eqref{4.6.1}
coincides with \eqref{4.5.1} owing to the
fact that  $r^{\alpha_{\cdot}
\beta_{\cdot}}_{t},\pi^{\alpha_{\cdot}
\beta_{\cdot}}_{t}$, and $w^{\alpha_{\cdot}
\beta_{\cdot}}_{t}$ are control adapted.
 It follows that (a.s.)
$$
x_{t}^{n}=x_{t}^{n}(\alpha_{\cdot})=
x^{nk}_{t}\quad\text{for all}\quad t\leq k/n,
$$
so that   (a.s.)
$$
\beta^{nk}_{t}(\alpha_{\cdot})=\beta(\alpha_{t},x^{n}_{k/n})
$$
for all $t\geq k/n$. Therefore, if $(k-1)/n\leq t<k/n$, then
$$
\beta^{n}_{t}(\alpha_{\cdot})=\beta^{n(k-1)}_{t}
(\alpha_{\cdot})=\beta(\alpha_{t},x^{n}_{(k-1)/n}),
$$

\begin{equation}
                                                         \label{1.19.1}
\beta^{n}_{s}:=\beta^{n}_{s}(\alpha_{\cdot})=\beta(\alpha_{s},
x^{n}_{\kappa_{n}(s)}),
\end{equation}
 and $x_{t}^{n}$
satisfies 
$$
x^{n}_{t}=x+\int_{0}^{t}
r^{n}_{s}
\sigma(\alpha_{s},\beta(\alpha_{s},
x^{n}_{\kappa_{n}(s)}),
 x^{n}_{s})\,
dw^{n          }_{s}
$$
\begin{equation}
                                                      \label{1.17.1}
+\int_{0}^{t}
[r^{n} _{s}]^{2}
\big[b(\alpha_{s},\beta(\alpha_{s},
x^{n}_{\kappa_{n}(s)}),x^{n}_{s}
)+\sigma(\alpha_{s},\beta(\alpha_{s},
x^{n}_{\kappa_{n}(s)}),
 x^{n}_{s})\pi^{n}_{s}
\big]\,
ds.
\end{equation}
with $(r,\pi,w)^{n}_{s}
=(r,\pi,w)^{\alpha_{\cdot}\beta^{n }_{\cdot} }_{s}$.  

 Introduce $\tau^{n}=\tau^{n}(\alpha_{\cdot})$ as the first exit time  
of $x^{n}_{t}=x^{n}_{t}(\alpha_{\cdot})$ from $G$ and set
$$
\phi^{n}_{t}=\phi^{\alpha_{\cdot}
\beta^{n}_{\cdot}x}_{t},\quad
\psi^{n}_{t}=\psi^{\alpha_{\cdot}
\beta^{n}_{\cdot}x}_{t} .
$$

 Observe that
by It\^o's formula
\begin{equation}
                                                          \label{2.10.1}
\hat{u}(x^{n}_{t\wedge\tau^{n}})
e^{-\phi^{n}_{t\wedge \tau^{n}}
-\psi^{n}_{t\wedge \tau^{n}}}
=\hat{u}(x)+\int_{0}^{t\wedge\tau^{n}}
[r^{n} _{s}]^{2}e^{-\phi^{n}_{s}-\psi^{n}_{s}}
 L^{\alpha_{s}\beta^{n}_{s}}
\hat{u}( x^{n}_{s})
\,ds+m^{n}_{t} ,  
\end{equation}
where $m^{n}_{s}$ is a martingale.
Here
  according to our assumptions on 
the uniform continuity in $x$ of the data 
and $D_{ij}\hat u(x)$ we  have that   
for $s<\tau^{n}$ (notice the change of $x^{n}_{s}$ to 
$x^{n}_{\kappa_{n}(s)}$)
$$
  L^{\alpha_{s}\beta^{n}_{s}}
\hat{u}(x^{n}_{s})\leq  
  a_{ij}(\alpha_{s},\beta(\alpha_{s},
x^{n}_{\kappa_{n}(s)}),x^{n}_{\kappa_{n}(s)})
D_{ij}\hat{u}(x^{n}_{\kappa_{n}(s)})
$$
$$
+  b_{i}(\alpha_{s},\beta(\alpha_{s},x^{n}_{\kappa_{n}(s)}),x^{n}_{\kappa_{n}(s)})
D_{i}\hat{u}(x^{n}_{\kappa_{n}(s)})
$$
$$
-  c(\alpha_{s},\beta(\alpha_{s},x^{n}_{\kappa_{n}(s)}),x^{n}_{\kappa_{n}(s)})
 \hat{u}(x^{n}_{\kappa_{n}(s)}) 
+\chi 
(x^{n}_{s}-x^{n}_{\kappa_{n}(s)}).
$$
where  $\chi (y)$ 
is a (nonrandom) bounded function on $\bR^{d}$ such that
$\chi (y)\to0$ as $y\to0$.
 All 
such functions will be denoted
by $\chi $ even if they may change from one occurrence to another.

Then  \eqref{1.17.3} shows that,
for $s<\tau^{n}$,
$$
  L^{\alpha_{s}\beta^{n}_{s}}
\hat{u}(x^{n}_{s})\leq \varepsilon+\chi 
(x^{n}_{s}-x^{n}_{\kappa_{n}(s)})
 - f(\alpha_{s},\beta(\alpha_{s}(x^{n}_{\kappa_{n}(s)}),
x^{n}_{\kappa_{n}(s)})
$$
$$
\leq \varepsilon+\chi 
(x^{n}_{s}-x^{n}_{\kappa_{n}(s)}) 
  -  f^{\alpha_{s}\beta^{n}_{s}}
( x^{n}_{s}),
$$
which along with \eqref{2.10.1}  
implies that, for
$$
\eta^{n\varepsilon}_{t}=\eta^{n\varepsilon}_{t}(\alpha_{\cdot})=
\delta_{1}^{-2}\int_{0}^{t\wedge\tau^{n}}
e^{-\phi^{n}_{s}-\psi^{n}_{s}}[
\varepsilon+
 \chi (x^{n}_{s}-x^{n}_{\kappa_{n}(s)})] \,ds ,
$$
\begin{equation}
                                           \label{1.17.5}
\kappa^{n\varepsilon}_{t\wedge\tau^{n}} 
 =\zeta^{n\varepsilon}_{t}+m^{n}_{t},
\end{equation}
where $\zeta^{n\varepsilon}_{t}$ is a decreasing process.

 Hence $\kappa^{n\varepsilon}_{t\wedge\tau^{n}}$ 
is at least a local supermartingale.  Owing to Lemmas
\ref{lemma 2.20.1} and \ref{lemma 2.10.1},
   \eqref{2.29.1} and
  \eqref{2.29.3}
hold.
It follows that the
 local supermartingale $\kappa^{n\varepsilon}_{t\wedge\tau^{n}}$  
 is, actually, a  supermartingale.

Furthermore,  Lemmas \ref{lemma 2.20.1}
 and \ref{lemma 2.10.1}, the boundedness
of $\chi$, its continuity, and the fact that $\chi(0)=0$
easily yield that 
\begin{equation}
                                           \label{1.17.6}
 \sup_{\alpha_{\cdot}\in\frA}E\int_{0}^{\tau^{n}(\alpha_{\cdot})}
e^{-\phi^{n}_{s}-\psi^{n}_{s}}
\chi (x^{n}_{s}(\alpha_{\cdot})-x^{n}_{\kappa_{n}(s)}
(\alpha_{\cdot})) \,ds\to0 
\end{equation}
as $n\to\infty$, which proves \eqref{3.3.1}. The lemma is proved.  

 For treating $\check{u}$ we use the following result.

\begin{lemma}  
                                                  \label{lemma 1.17.2}
Let $H[\check{u}]\geq0$  
  in $G$. Then for any $x\in\bR^{d}$, $\bbeta\in\bB$, and $\varepsilon>0$,
there exist a sequence $\alpha^{n}_{\cdot} 
\in\frA$, $n=1,2,...$, and a sequence of increasing continuous
$\{\cF_{t}\}$-adapted processes $\eta^{n\varepsilon}_{t}(\bbeta)$
 with $\eta^{n\varepsilon}_{0}(\bbeta)=0$
 such that 
the processes
$$
\kappa^{n\varepsilon}_{t} :=
\check{u}(x^{n}_{t})
e^{-\phi^{n}_{t}-\psi^{n}_{t}}
+\eta^{n\varepsilon}_{t} (\bbeta )
+\int_{0}^{t}
[r^{n}_{s}]^{2}f^{n}_{s}
( x^{n}_{s})
e^{-\phi^{n}_{s}-\psi^{n}_{s}}\,ds,
$$
where
\begin{equation}
                                               \label{5.6.4}
(x^{n}_{t},
\phi^{n}_{t},\psi^{n}_{t})
=(x_{t},\phi_{t},\psi_{t})
^{\alpha^{n}_{\cdot}\bbeta (\alpha^{n}_{\cdot})x},\quad
f^{n}_{t}( x)=
f^{\alpha^{n}_{t}\bbeta_{t}(\alpha^{n}_{\cdot}) }( x),
\quad r^{n}_{t}=
r^{\alpha^{n}_{\cdot}\bbeta (\alpha^{n}_{\cdot})}_{t},
\end{equation}
are submartingales on $[0,\tau^{\alpha^{n}_{\cdot}
\bbeta (\alpha^{n}_{\cdot})
  x}]$ and
\begin{equation}
                                                   \label{5.5.6}
\sup_{n}E 
\eta^{n\varepsilon}_{\infty}( \bbeta )<\infty,
\end{equation}
\begin{equation}
                                                   \label{5.5.7}
\nlimsup_{n\to\infty}
 E 
\eta^{n\varepsilon}_{\tau}(\bbeta)\leq
  N\varepsilon ,
\end{equation}
where $N$ 
is independent of $x$, $\bbeta$, and $\varepsilon$.
 
  Finally,
$$
\sup_{n}E 
\sup_{t\geq0}|
\kappa^{n\varepsilon}_{t\wedge\tau }|<\infty .
$$
\end{lemma}  

Proof. Owing to Assumptions
\ref{assumption 1.9.1}  
the function
$$
h(\alpha,x):=\inf_{\beta\in B}\big[   L^{\alpha\beta}\check{u}(x)+
  f^{\alpha\beta}(x)\big]
$$
is a finite Borel function of $x$ and is continuous with 
respect to $\alpha$.
Its $\sup$ over $A$ can be replaced with the $\sup$
over an appropriate countable subset of $A$ and since
$$
\sup_{\alpha\in A}h(\alpha,x)\geq0,
$$
similarly to how $\beta(\alpha,x)$ was defined
in the proof of Lemma \ref{lemma 1.16.1}, one can find 
 a Borel function $\bar{\alpha}(x)$ in such a way that
\begin{equation}
                                                \label{1.17.07}
\inf_{\beta\in B}\big[   L^{\bar{\alpha}(x)\beta}\check{u}(x)+
  f^{\bar{\alpha}(x)\beta}(x)\big]\geq  -\varepsilon 
\end{equation}
in $G$. If $x\not\in G$ we set $\bar{\alpha}(x)=\alpha^{*}$,
where $\alpha^{*}$ is a fixed element of $A$.  

After that we need some processes which we introduce recursively.
Fix $x$ and set $\alpha^{n0}_{t}\equiv\bar{\alpha}(x )$. Then define
$x^{n0}_{t}$, $t\geq0$, as a unique solution of the equation 
$$
x _{t}=x+\int_{0}^{t}r_{s}^{\alpha^{n0}_{\cdot}
\bbeta  (\alpha^{n0}_{\cdot })}
\sigma(\alpha^{n0}_{s},
\bbeta _{s}(\alpha^{n0}_{\cdot }),x _{s})\,dw_{s}^{\alpha^{n0}_{\cdot}
\bbeta  (\alpha^{n0}_{\cdot })} 
$$
$$
+\int_{0}^{t}[r_{s}^{\alpha^{n0}_{\cdot}
\bbeta  (\alpha^{n0}_{\cdot })}]^{2}\big[
b(\alpha^{n0}_{s},
\bbeta _{s}(\alpha^{n0}_{\cdot }) ,x _{s})
+\sigma(\alpha^{n0}_{s},
\bbeta _{s}(\alpha^{n0}_{\cdot }),x _{s})
\pi_{s}^{\alpha^{n0}_{\cdot}
\bbeta  (\alpha^{n0}_{\cdot })}\big]\,ds.
$$
For $k\geq1$ introduce $\alpha^{nk}_{t}$ so that
$$
\alpha^{nk}_{t}=\alpha^{n(k-1)}_{t}\quad\text{for}\quad t<k/n,
$$
$$
\alpha^{nk}_{t}=\bar{\alpha}(x^{n(k-1)}_{k/n})
\quad\text{for}\quad t\geq k/n,
$$
where $x^{n(k-1)}_{t}$ is a unique solution of
 $$
x _{t}=x+\int_{0}^{t}r_{s}^{\alpha^{n(k-1)}_{\cdot}
\bbeta  (\alpha^{n(k-1)}_{\cdot })}
\sigma(\alpha^{n(k-1)}_{s},
\bbeta _{s}(\alpha^{n(k-1)}_{\cdot }), x _{s})\,dw_{s}^{\alpha^{n(k-1)}_{\cdot}
\bbeta  (\alpha^{n(k-1)}_{\cdot })} 
$$
$$
+\int_{0}^{t}[r_{s}^{\alpha^{n(k-1)}_{\cdot}
\bbeta  (\alpha^{n(k-1)}_{\cdot })}]^{2}\big[b(\alpha^{n(k-1)}_{s},
\bbeta _{s}(\alpha^{n(k-1)}_{\cdot }) ,x _{s})
$$
\begin{equation}
                                                      \label{5.5.1}
+\sigma(\alpha^{n(k-1)}_{s},
\bbeta _{s}(\alpha^{n(k-1)}_{\cdot }), x _{s})
\pi_{s}^{\alpha^{n(k-1)}_{\cdot}
\bbeta  (\alpha^{n(k-1)}_{\cdot })}\big]\,ds.
\end{equation}

As in the proof of Lemma \ref{lemma 1.16.1} one can 
show that the above definitions
make sense as well as the definition
\begin{equation}
                                                \label{5.5.2}
\alpha^{n}_{t}=\alpha^{n(k-1)}_{t}\quad\text{for}\quad t<k/n.
\end{equation}
Next, by definition $x^{n}_{t}=
 x_{t}^{\alpha_{\cdot}^{n}\bbeta(\alpha^{n}_{\cdot})x}$ satisfies
$$
x_{t}=x+\int_{0}^{t}r_{s}^{\alpha^{n}_{\cdot} \bbeta (
\alpha^{n}_{\cdot})}
\sigma(\alpha^{n}_{s},\bbeta_{s}(
\alpha^{n}_{\cdot}),x_{s})\,dw_{s}^{\alpha^{n}_{\cdot} \bbeta (
\alpha^{n}_{\cdot})}
$$
$$
+\int_{0}^{t}[r_{s}^{\alpha^{n}_{\cdot} \bbeta (
\alpha^{n}_{\cdot})}]^{2}\big[b(\alpha^{n}_{s},\bbeta_{s}(
\alpha^{n}_{\cdot}) ,x_{s})
+\sigma(\alpha^{n}_{s},\bbeta_{s}(
\alpha^{n}_{\cdot}),x_{s})\pi_{s}^{\alpha^{n}_{\cdot} \bbeta (
\alpha^{n}_{\cdot})}\big]\,ds.
$$
Equation  \eqref{5.5.2}   and the definitions
of $\bB$ and of control adapted processes show that
$x^{n}_{t}$ satisfies  \eqref{5.5.1} for $t\leq k/n$.
Hence, (a.s.) $x^{n}_{t}=x^{n(k-1)}_{t}$ for all $t\leq k/n$
and (a.s.) for all $t\geq0$, $\alpha^{n}_{t}
=\bar{\alpha}(x^{n}_{\kappa_{n}(t)})$ and 
$$
x^{n}_{t}=x+\int_{0}^{t}r^{n}_{s}\sigma(\bar{\alpha}(x^{n}_{\kappa_{n}(s)}),
\bbeta _{s}(\alpha^{n}_{\cdot }),x^{n}_{s})\,dw^{n}_{s}
$$
$$
+\int_{0}^{t}[r^{n}_{s}]^{2}\big[b(\bar{\alpha}(x^{n}_{\kappa_{n}(s)}),
\bbeta _{s}(\alpha^{n}_{\cdot }), x^{n}_{s})
+\sigma(\bar{\alpha}(x^{n}_{\kappa_{n}(s)}),
\bbeta _{s}(\alpha^{n}_{\cdot }),x^{n}_{s})\pi^{n}_{s}\big]\,ds,
$$
where $(r,\pi,w)^{n}_{s} =(r,\pi,w)_{s}^{\alpha^{n}_{\cdot} \bbeta (
\alpha^{n}_{\cdot})}$.

Now, 
introduce $\tau^{n}$ as the first exit time
of $x^{n}_{t}$ from $G$, set
$$
\beta^{n}_{s}=\bbeta_{s}(\alpha^{n}_{\cdot}),\quad
\phi^{n}_{t}=\phi^{\alpha^{n}_{\cdot}
\beta^{n}_{\cdot} x}_{t}
,\quad
\psi^{n}_{t}=\psi^{\alpha^{n}_{\cdot}
\beta^{n}_{\cdot} x}_{t},\quad r^{n}_{s}=
r^{\alpha^{n}_{\cdot}
\beta^{n}_{\cdot} } ,
$$
 and observe that
by It\^o's formula  

$$
\check{u}(x^{n}_{t\wedge\tau^{n}})
e^{-\phi^{n}_{t\wedge \tau^{n}}
-\psi^{n}_{t\wedge \tau^{n}}}
=\check{u}(x)+\int_{0}^{t\wedge\tau}
[r^{n}_{s}]^{2}e^{-\phi^{n}_{s}-\psi^{n}_{s}}
  L^{\alpha^{n}_{s} 
\beta^{n}_{s}}
\check{u}(x^{n}_{s})
\,ds+m^{n}_{t},
$$
where $m^{n}_{s}$ is a martingale  and,
for $s<\tau^{n}$,
$$
  L^{\alpha^{n}_{s}
\beta^{n}_{s}}
\check{u}(x^{n}_{s})=  
a_{ij}(\bar{\alpha}(x^{n}_{\kappa_{n}(s)}),\beta^{n}_{s}
,x^{n}_{s})
D_{ij}\check{u}(x^{n}_{s})
$$
$$
+  b_{i}(\bar{\alpha}(x^{n}_{\kappa_{n}(s)}),\beta^{n}_{s}
,x^{n}_{s})
D_{i}\check{u}(x^{n}_{s})
-  c(\bar{\alpha}(x^{n}_{\kappa_{n}(s)}),\beta^{n}_{s}
,x^{n}_{s})
 \check{u}(x^{n}_{s}).
$$
Similarly to the proof of Lemma \ref{lemma 1.16.1} we derive from 
\eqref{1.17.07} that,
for $s<\tau^{n}$,
$$
\bar L^{\alpha^{n}_{s}\beta^{n}_{s}}
\check{u}(x^{n}_{s})\geq -\varepsilon-\chi(x^{n}_{s}-x^{n}_{\kappa_{n}(s)})  -
f(\bar{\alpha}(x^{n}_{\kappa_{n}(s)}),\beta^{n}_{s},x^{n}_{\kappa_{n}(s)})
$$
$$
= -\varepsilon-\chi(x^{n}_{s}-x^{n}_{\kappa_{n}(s)})  
 -  f^{\alpha^{n}_{s}\beta^{n}_{s}}
( x^{n}_{s}),
$$
where   $\chi (y)$ 
are   (nonrandom) bounded functions on $\bR^{d}$ such that
$\chi (y)\to0$ as $y\to0$.
It follows that
\begin{equation}
                                           \label{1.17.8}
\check{u}(x^{n}_{t\wedge\tau^{n}})
e^{-\phi^{n}_{t\wedge\tau^{n}}
-\psi^{n}_{t\wedge\tau^{n}}}
+\int_{0}^{t\wedge\tau^{n}}[r^{n}_{s}]^{2}
f^{\alpha^{n}_{s}\beta^{n}_{s}   }
( x^{n}_{s})
e^{-\phi^{n}_{s}-\psi^{n}_{s}}\,ds+\eta^{n}_{t}
 =\zeta_{t}+m^{n}_{t},
\end{equation}
where $\zeta_{t}$ is an  increasing process and
$$
\eta^{n}_{t}=\eta^{n}_{t}(\bbeta)=
\delta_{1}^{-2}\int_{0}^{t\wedge\tau^{n}}
e^{-\phi^{n}_{s}-\psi^{n}_{s}}[\varepsilon
 +
 \chi_{\varepsilon}(x^{n}_{s}-x^{n}_{\kappa_{n}(s)})] \,ds.
$$
 Hence the left-hand side of \eqref{1.17.8}
is a local submartingale and we finish the proof
in the same way as the proof of Lemma \ref{lemma 1.16.1}.
 The lemma is proved.

{\bf Proof of Theorem \ref{theorem 2.11.1}}. (i)
First we fix $x\in\bR^{d}$,
$\alpha_{\cdot}\in\frA$, and $\varepsilon>0$,
 take $\beta^{n}_{\cdot}(\alpha_{\cdot})$
form Lemma \ref{lemma 1.16.1} and prove that  the $\frB$-valued
functions defined on $\frA$ by
 $\bbeta^{n}(\alpha_{\cdot})=\beta^{n}_{\cdot}(\alpha_{\cdot})$
belong to $\bB$. To do that observe that if \eqref{4.5.4} holds and $T\leq 1/n$, then (a.s.)
$\beta^{n0}_{t}(\alpha^{1}_{\cdot})=\beta^{n0}_{t}(\alpha^{2}_{\cdot})$
for almost all $t\leq T$. By definition also (a.s.)
$$
(r,\pi,w)^{\alpha^{1}_{\cdot}\beta_{\cdot}^{n0}(\alpha^{1}_{\cdot})}_{s}
=(r,\pi,w)^{\alpha^{2}_{\cdot}\beta_{\cdot}^{n0}(\alpha^{2}_{\cdot})}_{s}
\quad\text{for almost all}\quad s\leq T.
$$
By uniqueness of solutions of \eqref{5.11.1} (see Assumption
\ref{assumption 5.23.2}),
 the processes $x_{t}^{n1}$ found from \eqref{4.5.1}
for $\alpha_{\cdot}=\alpha^{1}_{\cdot}$
and for $\alpha_{\cdot}=\alpha^{2}_{\cdot}$ coincide (a.s.)
for all $t\leq T$.

If \eqref{4.5.4} holds and $1/n<T\leq 2/n$, then by the above
solutions of \eqref{4.5.1} for $\alpha_{\cdot}=\alpha^{1}_{\cdot}$
and for $\alpha_{\cdot}=\alpha^{2}_{\cdot}$ coincide (a.s.)
for $t=1/n$ and then (a.s.)
$\beta^{n1}_{t}(\alpha^{1}_{\cdot} )=
\beta^{n1}_{t}(\alpha^{2}_{\cdot} )$ not only for all $t<1/n$
but also
for all $t\geq1/n$, 
which implies that (a.s.)
$$
(r,\pi,w)^{\alpha^{1}_{\cdot}\beta_{\cdot}^{n1}(\alpha^{1}_{\cdot})}_{s}
=(r,\pi,w)^{\alpha^{2}_{\cdot}\beta_{\cdot}^{n1}(\alpha^{2}_{\cdot})}_{s}
\quad\text{for almost all}\quad s\leq T
$$
and again the processes $x_{t}^{n }$ found from \eqref{4.5.1}
for $\alpha_{\cdot}=\alpha^{1}_{\cdot}$
and for $\alpha_{\cdot}=\alpha^{2}_{\cdot}$ coincide (a.s.)
for all $t\leq T$.

By induction we get that if \eqref{4.5.4} holds for a $T\in(0,\infty)$
and we define
$k $ as the integer such that $k/n<T\leq(k+1)/n$, then (a.s.)
\begin{equation}
                                               \label{4.5.5}
\beta^{n}_{t}(\alpha^{1}_{\cdot})=\beta^{nk}_{t}(\alpha^{1}_{\cdot} )=
\beta^{nk}_{t}(\alpha^{2}_{\cdot} )=
\beta^{n}_{t}(\alpha^{2}_{\cdot})\quad\text{for almost all}\quad
t< (k+1)/n ,
\end{equation}
$$
(r,\pi,w)^{\alpha^{1}_{\cdot}\beta_{\cdot}^{nk}(\alpha^{1}_{\cdot})}_{s}
=(r,\pi,w)^{\alpha^{2}_{\cdot}\beta_{\cdot}^{nk}(\alpha^{2}_{\cdot})}_{s}
\quad\text{for almost all}\quad s\leq T
$$
and  the processes $x_{t}^{n }$ found from \eqref{4.5.1}
for $\alpha_{\cdot}=\alpha^{1}_{\cdot}$
and for $\alpha_{\cdot}=\alpha^{2}_{\cdot}$ coincide (a.s.)
for all $t\leq T$.
This means that $\bbeta^{n} \in\bB$ indeed.

Furthermore, by the supermartingale property of 
$\kappa^{n\varepsilon}_{t}(\alpha_{\cdot})$,
  we have
$$
\hat{u}(x)\geq
E_{x}^{\alpha_{\cdot}\bbeta^{n}(\alpha_{\cdot})}\big[
g(x_{\tau})e^{-\phi_{\tau}-\psi_{\tau} }
+\int_{0}^{\tau}r^{2}_{t}
 f( x_{t})  e^{-\phi_{t}-\psi_{t} }\,dt\big]
-E 
\eta^{n\varepsilon}_{\tau}
(\alpha_{\cdot})  ,
$$
which owing to \eqref{3.3.1} yields
$$
\hat{u}(x)\geq\nliminf_{n\to\infty}\sup_{\alpha_{\cdot}\in\frA}
E_{x}^{\alpha_{\cdot}\bbeta^{n} (\alpha_{\cdot})}
\big[
 \int_{0}^{\tau}
r^{2}_{t}f(x_{t}) e^{-\phi_{t}-\psi_{t}}\,dt
+g(x_{\tau})e^{-\phi_{\tau}-\psi_{\tau}} \big]
-N\varepsilon .
$$
 
In light of the arbitrariness of $\varepsilon$ we conclude
$\hat u\geq v$
  and assertion (i) is proved.

(ii)
 Similarly to the above argument,
for any $\bbeta\in\bB$,
$$
\check{u}(x)\leq E_{x}^{\alpha^{n}_{\cdot}\bbeta (\alpha^{n}_{\cdot})}
\big[\int_{0}^{\tau}r^{2}_{t}f( x_{t})
e^{-\phi_{t}-\psi_{t} }\,dt
+g(x_{\tau})e^{-\phi_{\tau}-\psi_{\tau}}\big]
+E\eta^{n\varepsilon}_{\tau}(\bbeta ) .
$$

It follows that
$$
\check{u}(x)\leq \sup_{\alpha_{\cdot}\in\frA}
E_{x}^{\alpha _{\cdot}\bbeta(\alpha _{\cdot})}
\big[\int_{0}^{\tau}r^{2}_{t}f( x_{t})
e^{-\phi_{t} -\psi_{t} }\,dt
+g(x_{\tau})e^{-\phi_{\tau}-\psi_{\tau} }\big]
+\nlimsup_{n\to\infty}
E 
\eta^{n\varepsilon}_{\tau}
(\bbeta)   
$$
$$
 \leq \sup_{\alpha_{\cdot}\in\frA}
E_{x}^{\alpha _{\cdot}\bbeta(\alpha _{\cdot})}
\big[\int_{0}^{\tau}r^{2}_{t}f( x_{t})
e^{-\phi_{t}-\psi_{t}}\,dt
+g(x_{\tau})e^{-\phi_{\tau}-\psi_{\tau}}\big]+N\varepsilon,
$$
which in light of the arbitrariness of $\varepsilon$ and $\bbeta\in\bB$
 finally yields that $\check u\leq v$.

This proves assertion (ii). Assertion (iii) is an obvious
consequence of (i) and (ii).
The theorem is proved.

\mysection{The case of uniformly nondegenerate processes}
                                                         \label{section 4.6.1}

As in Section \ref{section 2.11.1}  we   replace Assumption 
\ref{assumption 1.9.1} (iii)
with Assumptions \ref{assumption 2.16.1}
and \ref{assumption 5.23.2}, and we  assume that we are given two  functions
 $\hat{u},\check{u} \in W^{2}_{d,loc}(G)\cap C(\bar{G})$.

 In that case we have the  following.
\begin{theorem}
                                         \label{theorem 2.16.1}
 (i) If $H[\hat{u}]\leq0$ (a.e.)
  in $G$ 
and $\hat{u}\geq g$ on $\partial G$, then $v\leq\hat{u}$ in $\bar{G}$.

(ii) If $H[\check{u}]\geq0$ (a.e.)
  in $G$ 
and $\check{u}\leq g$  on $\partial G$,
 then $v\geq\check{u}$ in $\bar{G}$.

(iii) If $\hat{u}$ and $\check{u}$ are as in (i) and (ii) and
  $\hat{u}=\check{u}$, then
$v=\hat{u}$ and $v$ is independent of the choice of the probability space,
filtration, $r,\pi$, and $w$. 

\end{theorem}

Proof. (i) We basically repeat the proof of Theorem 4.1
 of \cite{Kr_13_1} with considerable
simplifications made possible due to our
assumptions. It is well known that there exists
a sequence $\hat{u}_{n}\in C^{2}(\bar{G})$ such that
$\hat{u}_{n}\to\hat{u}$ in $C(\bar{D})$ and in $W^{2}_{d}(G')$
for any subdomain $G'\subset \bar{G}'\subset G$.
Introduce 
$\hat{h}_{n}=H[\hat{u}_{n}]$,
$$
f^{\alpha\beta}_{n}(x)= f^{\alpha\beta} (x)
 - \hat{h}_{n}(x),  
$$
and observe that owing to our continuity assumptions 
on $\sigma,b,c,f$, the functions $\hat{h}_{n}$ and 
$f^{\alpha\beta}_{n}(x)$ are continuous in $x$ uniformly with respect to
$\alpha,\beta$ and
$$
\supinf_{\alpha\in A\,\,\beta\in B}
\big[
L^{\alpha\beta}_{n}\hat{u}_{n}( x)+
f^{\alpha\beta}_{n}( x)\big]=0
$$
in $G$. By Theorem \ref{theorem 2.11.1},
for any subdomain $G_{1}\subset \bar{G}_{1}\subset G$
 we have in  $G_{1}$ that
\begin{equation}
                                                         \label{2.16.3}
\hat{u}_{n}(x)\geq\infsup_{\bbeta\in\bB\,\,\alpha_{\cdot}\in\frA}
E_{x}^{\alpha_{\cdot}\bbeta(\alpha_{\cdot})}\big[
\hat{u}_{n}(x_{\tau_{1}})e^{-\phi_{\tau_{1}}
-\psi_{\tau_{1}}}
+\int_{0}^{\tau_{1}}r^{2}_{t}
 f_{n}( x_{t}) e^{-\phi_{t}-\psi_{t}}\,dt \big],
\end{equation}
where
$$
\tau_{1}^{\alpha_{\cdot}\beta_{\cdot}x}
=\inf\{t\geq0:x_{t}^{\alpha_{\cdot}\beta_{\cdot}x}\not\in G_{1}\}.
$$
 
 Notice that
$$
E^{\alpha_{\cdot}\beta_{\cdot}}_{ x}
 |\hat{u}_{n}(x_{\tau_{1}})- \hat{u} (x_{\tau_{1}})|
e^{-\phi_{\tau_{1}}
-\psi_{\tau_{1}}} \leq \sup_{G}|\hat{u}_{n}-\hat{u}|.
$$
While estimating
$$
I_{n}(x):=E^{\alpha_{\cdot}\beta_{\cdot}}_{ x}\int_{0}^{\tau_{1}}
|f_{n}-f|(x_{t})e^{-\phi_{t}-\psi_{t}}\,dt
=E^{\alpha_{\cdot}\beta_{\cdot}}_{ x}e^{-\psi_{\tau_{1}}}\int_{0}^{\tau_{1}}
|f_{n}-f|(x_{t})e^{-\phi_{t}}\,dt.
$$
Girsanow's theorem allows us to concentrate on $\pi\equiv0$ and then
  the Alexandrov estimate guarantees that
$$
I_{n}(x)\leq N\|\hat{h}_{n}\|_{\cL_{d}(G_{1})}
=N\| H[\hat{u}_{n}]-H[\hat{u} ]\|_{\cL_{d}(G_{1})}
\leq N\|\hat{u}_{n} - \hat{u}\|_{W^{2}_{d}(G_{1})},
$$
where the constants $N$ are independent of $n$ (and $x$).

Hence by letting $n\to\infty$ in \eqref{2.16.3} we obtain
that for $k=1$
\begin{equation}
                                                         \label{2.16.4}
\hat{u} (x)\geq\infsup_{\bbeta\in\bB\,\,\alpha_{\cdot}\in\frA}
E_{x}^{\alpha_{\cdot}\bbeta(\alpha_{\cdot})}\big[
\hat{u} (x_{\tau_{k}})e^{-\phi_{\tau_{k}} 
-\psi_{\tau_{k}}}
+\int_{0}^{\tau_{k}}
 f ( x_{t}) e^{-\phi_{t}-\psi_{t}}\,dt \big],
\end{equation}
Where $\tau_{k}^{\alpha_{\cdot}\beta_{\cdot}x}$ are defined
as the first exit times of the processes
$x_{t}^{\alpha_{\cdot}\beta_{\cdot}x}$ from an expanding
sequence of subdomains $G_{k}\subset\bar{G}_{k}\subset G$
such that $\cup_{k}G_{k}= G$.

 By letting $k\to\infty$ in \eqref{2.16.4} and 
repeating the proof of Theorem 2.2 of \cite{Kr_13_1}
given there in Section 6 we get that $\hat{u}\geq v$ in $G$
as stated.
Observe that in our situation in the proof
of Theorem 2.2 of \cite{Kr_13_1} we need not mollify $f^{\alpha\beta}(x)$
because by assumption it is uniformly continuous in $x$.

The proof of assertion (ii) is quite similar and 
as usual assertion (iii) is obtained by  simply
combining assertions (i) and (ii). The theorem is proved.

\mysection{A general approximation result from above}
                                              \label{section 2.25.1} 

In this section
we suppose that all assumptions in Section \ref{section 2.26.3} are
satisfied.
Set 
$$
A_{1}=A 
$$
 and let $A_{2}$ be a 
separable metric space having no common points with $A_{1}$.

\begin{assumption}
                                         \label{assumption 4.29.2}
The functions $ 
\sigma^{\alpha\beta}( x)$,
$ b^{\alpha\beta}( x)$, $ 
c^{\alpha\beta}( x)$, and
$ f^{\alpha\beta}( x)$ 
 are also defined on
  $A_{2}\times B \times\bR^{d}$ in such a way that they are
{\em independent\/} of $\beta$ (on
  $A_{2}\times B \times\bR^{d}$) and 
the assumptions in Section \ref{section 2.26.3}  are satisfied,
  of course, with $A_{2}$ in place of $A$.

\end{assumption}

Define
$$
\hat{A}=A_{1}\cup A_{2}.
$$

Then we introduce $\hat{\frA}$ as the set of progressively measurable
$\hat{A}$-valued processes and $\hat{\bB}$ as the set
of $\frB$-valued functions $ \bbeta(\alpha_{\cdot})$
on $\hat{\frA}$ such that,
for any $T\in[0,\infty)$ and any $\alpha^{1}_{\cdot},
\alpha^{2}_{\cdot}\in\hat{\frA}$ satisfying
$$
P(  \alpha^{1}_{t}=\alpha^{2}_{t} 
 \quad\text{for almost all}\quad t\leq T)=1,
$$
we have
$$
P(  
 \bbeta_{t}(\alpha^{1}_{\cdot})=\bbeta_{t}(\alpha^{2}_{\cdot}) 
 \quad\text{for almost all}\quad t\leq T)=1.
$$

We fix an element $\alpha^{*}\in A_{1}$ and for
$\alpha_{\cdot}\in\hat{\frA}$ define
$$
(p\alpha)_{t}=\alpha_{t}\quad\text{if}
\quad\alpha_{t}\in A_{1},\quad
(p\alpha)_{t}=\alpha^{*}\quad\text{if}
\quad\alpha_{t}\in A_{2}.
$$
By using this projection operator we extend
 $(w,r,\pi)^{\alpha_{\cdot}\beta_{\cdot}}_{t}$
originally defined for 
$\alpha_{\cdot}\in\frA$ and $\beta_{\cdot}\in\frB$
as
\begin{equation}
                                                    \label{2.22.1}
(w,r,\pi)^{\alpha_{\cdot}\beta_{\cdot}}_{t}
=(w,r,\pi)^{p\alpha_{\cdot}\beta_{\cdot}}_{t}
\end{equation}
thereby now defined for
$\alpha_{\cdot}\in\hat{\frA}$ and $\beta_{\cdot}\in\frB$.

Next, take a constant $K\geq0$ and set 
$$
v_{K}(x)=\infsup_{\bbeta\in\hat{\bB}\,\,\alpha_{\cdot}\in\hat{\frA}}
v^{\alpha_{\cdot}\bbeta(\alpha_{\cdot})}_{K}(x),
$$
where
$$
v^{\alpha_{\cdot}\beta _{\cdot} }_{K}(x)=
E_{x}^{\alpha_{\cdot}\beta _{\cdot} }\big[\int_{0}^{\tau}
r^{2}_{t} f_{K} ( x_{t})e^{-\phi_{t}-\psi_{t}}\,dt
+g(x_{\tau})e^{-\phi_{\tau}-\psi_{\tau}
 }\big]
$$
$$
= :v^{\alpha_{\cdot}\beta _{\cdot} }  (x)-K
E_{x}^{\alpha_{\cdot}\beta _{\cdot} }\int_{0}^{\tau}
r^{2}_{t}I_{\alpha_{t}\in A_{2}}e^{-\phi_{t} 
-\psi_{t}}\,dt,
$$
$$
 f^{\alpha\beta}_{K}( x)=f^{\alpha\beta}( x)-KI_{\alpha\in A_{2}}.
$$
Notice that, obviously,
$$
v (x)=\infsup_{\bbeta\in\hat{\bB}\,\,\alpha_{\cdot}\in\hat{\frA}}
v^{p\alpha_{\cdot}\bbeta(p\alpha_{\cdot})} (x).
$$

These definitions make sense  
owing to Remark \ref{remark 2.20.1}, which also implies that
  $v^{\alpha_{\cdot}\beta _{\cdot} }_{K}$
and $v^{\alpha_{\cdot}\beta _{\cdot} }$  
and bounded in $\bar{G}$.

\begin{theorem}
                                             \label{theorem 1.20.1}
We have
 $v_{K}\to v$ uniformly on $\bar{G}$ as $K\to\infty$.
\end{theorem}

\begin{lemma}
                                    \label{lemma 3.24.1}
Assume that $\pi\equiv0$. Then
there exists a constant $N$
such that for any
$\alpha_{\cdot}\in\hat{\frA}$, $\beta_{\cdot}\in\frB$,
 $x\in\bR^{d}$,  $T\in[0,\infty)$, and
stopping time $\gamma$
$$
E^{\alpha_{\cdot}\beta_{\cdot}}_{x}
\sup_{t\leq T\wedge\gamma}|x_{t}-y_{t}|\leq
Ne^{NT} \big(E^{\alpha_{\cdot}\beta_{\cdot}}_{x}\int_{0}^{
T\wedge\gamma}I_{\alpha_{t}\in A_{2}}\,dt\big)^{1/2} ,
$$
where 
$$
y^{\alpha_{\cdot} \beta_{\cdot}x}_{t}
=x^{p\alpha_{\cdot} \beta_{\cdot}x}_{t}.
$$
\end{lemma}

Proof. For simplicity of notation we drop the
superscripts $\alpha_{\cdot} ,\beta_{\cdot},x$.
Observe that $x_{t}$ and $y_{t}$ satisfy
$$
x_{t}=x+\int_{0}^{t}r_{s}\sigma^{\alpha_{s}\beta_{s}}( x_{s})
\,dw_{s}+\int_{0}^{t}r^{2}_{s}
b^{\alpha_{s}\beta_{s}}( x_{s})\,ds,
$$
$$
y_{t}=x+\int_{0}^{t}r_{s}\sigma^{\alpha_{s}\beta_{s}}( y_{s})
\,dw_{s}+\int_{0}^{t}r_{s}^{2}b^{\alpha_{s}\beta_{s}}( y_{s})
\,ds+\eta_{t},
$$
where $\eta_{t}=I_{t}+J_{t}$,
$$
I_{t}=\int_{0}^{t}r_{s}
[\sigma^{p\alpha_{s}\beta_{s}}( y_{s})
-\sigma^{\alpha_{s}\beta_{s}}( y_{s})]\,dw_{s},
$$
$$
J_{t}=\int_{0}^{t}r^{2}_{s}[b^{p\alpha_{s}\beta_{s}}( y_{s})
-b^{\alpha_{s}\beta_{s}}( y_{s})]\,ds.
$$

By Theorem II.5.9 of \cite{Kr77} (where we replace
the processes $x_{t}$ and $\tilde{x}_{t}$ with appropriately stopped ones)
for any $T\in[0,\infty)$ and any stopping time $\gamma
 $
\begin{equation}
                                                       \label{1.22.04}
E\sup_{t\leq T\wedge \gamma}|x _{t}-y _{t}|^{2}\leq
Ne^{NT}E\sup_{t\leq T\wedge \gamma}|\eta _{t}|^{2},
\end{equation}
where $N$ depends only on $K_{1}$  and $d$, 
which by Theorem III.6.8
of \cite{Kr95} leads to
\begin{equation}
                                                       \label{1.23.01}
E\sup_{t\leq T\wedge\gamma }|x _{t}-y _{t}|
\leq Ne^{NT}E\sup_{t\leq T\wedge\gamma }|\eta_{t}|
\end{equation}
with the constant $N$ being three times the one from \eqref{1.22.04}.

By using Davis's inequality   we see that for any $T\in[0,\infty)$
$$
E\sup_{t\leq T\wedge\gamma }|I _{t}| 
\leq NE\big(\int_{0}^{T\wedge\gamma}
I_{\alpha _{s}\in A_{2}}\,ds\big)^{1/2}
\leq N\big(E\int_{0}^{T\wedge\gamma}
I_{\alpha _{s}\in A_{2}}\,ds\big)^{1/2}.
$$

Furthermore, almost obviously
$$
E\sup_{t\leq T\wedge\gamma }|J_{t}| 
\leq N E\int_{0}^{T\wedge\gamma}
I_{\alpha _{s}\in A_{2}}\,ds 
\leq NT^{1/2}\big(E\int_{0}^{T\wedge\gamma}
I_{\alpha _{s}\in A_{2}}\,ds\big)^{1/2}
$$
and this in combination with \eqref{1.23.01} proves the lemma.

{\bf Proof of Theorem \ref{theorem 1.20.1}}.
Without losing generality we may
assume  that $g\in C^{3}(\bR^{d})$ since
the functions of this class   uniformly  
approximate in $\bar{G}$  any $g$ which is   continuous 
in $\bR^{d}$. Then notice that by It\^o's formula
for $g\in C^{3}(\bR^{d})$  
we have
$$
E_{x}^{\alpha_{\cdot}\beta _{\cdot} }\big[\int_{0}^{\tau}
r_{t}^{2} f_{K} ( x_{t})e^{-\phi_{t}-\psi_{t} }\,dt
+g(x_{\tau})e^{-\phi_{\tau}-\psi_{\tau}
 }\big]
$$
$$
=g(x)+E_{x }^{\alpha _{\cdot} \beta 
 _{\cdot} } \int_{0}^{\tau}r_{t}^{2}[ 
 \hat{f} ( x_{t} )-KI_{\alpha_{t}\in A_{2}}] 
e^{-\phi_{t}-\psi_{t}}\,dt,
$$
where 
$$
\hat{f} ^{\alpha\beta}( x ):=f^{\alpha\beta} ( x ) 
+ L^{\alpha\beta}g( x),
$$
which
is bounded and, for $(\alpha,\beta)\in\hat A\times B$, is
 uniformly continuous in $x$ uniformly with respect to $\alpha,\beta$.
This argument   shows that without losing generality
we may (and will) also assume that $g=0$.

Next, since   $\frA\subset\hat{\frA}$ and for $\alpha_{\cdot}\in
\hat{\frA}$
and $\bbeta \in\hat{\bB}$ we have $\bbeta(\alpha_{\cdot})\in\frB$,
it holds that
 $$ 
v_{K}\geq v.
$$

To estimate $v_{K}$ from above, take $\bbeta\in\bB$ and define 
$\hat{\bbeta} \in\hat{\bB}$ by
\begin{equation}
                                                              \label{2.22.2}
\hat{\bbeta}_{t}(\alpha_{\cdot})=\bbeta_{t}(p\alpha_{\cdot}).
\end{equation}
Also take any sequence $x^{n}
 \in \bar{G}$, $n=1,2,...$, recall that
$r_{t}^{\alpha_{\cdot}\beta_{\cdot}}\geq \delta_{1}$, and find a sequence
$\alpha^{n}_{\cdot}\in\hat{\frA}$ such that
$$
v_{K}(x^{n}) \leq  \sup_{\alpha\in\hat{\frA}}
E_{x^{n}}^{\alpha_{\cdot}\hat{\bbeta}(\alpha_{\cdot})}
 \int_{0}^{\tau}r^{2}_{t}
 f _{K} ( x_{t})e^{-\phi_{t}-\psi_{t}}\,dt
$$
\begin{equation}
                                                    \label{2.22.6}
=1/n+v^{\alpha^{n}_{\cdot}\hat{\bbeta}(\alpha^{n}_{\cdot})}(x^{n})
-K\delta_{1}^{2}E\int_{0}^{\tau^{n} }I_{\alpha^{n}_{t}\in A_{2}}
e^{-\phi^{n}_{t}-\psi^{n}_{t}}\,dt,
\end{equation}
where
$$
(\tau^{n},\phi^{n}_{t},\psi^{n}_{t} )
=( \tau,\phi _{t},\psi_{t} )
^{\alpha^{n}_{\cdot}\hat{\bbeta}(\alpha^{n}_{\cdot})x^{n}}.
$$
It follows   that
 there is a constant $N$
independent of $n$ and $K$ such that
\begin{equation}
                                                             \label{2.22.4}
E\int_{0}^{\tau^{n} }I_{\alpha^{n}_{t}\in A_{2}}
e^{-\phi^{n}_{t}-\psi^{n}_{t}}\,dt\leq N/K.
\end{equation}
Below by $N$ we denote generic constants independent of $n$
and $K$ (and $T$ once it appears).

We want to estimate the difference
\begin{equation}
                                                             \label{2.22.3}
v^{\alpha^{n}_{\cdot}\hat{\bbeta}(\alpha^{n}_{\cdot})}(x^{n})
-v^{p\alpha_{\cdot}^{n}\bbeta(p\alpha_{\cdot}^{n})}(x^{n}).
\end{equation}
Observe that 
in the expression of this difference by the definition
through the mathematical expectations of certain
quantities  the processes
$\psi_{t}$ involved are just the same, thanks to \eqref{2.22.1}
and \eqref{2.22.2}. This allows us to rewrite
the mathematical expectations similarly to how it is done
in Remark \ref{remark 2.20.1} and then by using Girsanov's theorem
  allows us to assume that $\pi_{t}\equiv0$, at the expense
that the underlying probability measures will now depend on $n$.
However, for simplicity of notation we keep the symbol $E$
for expectations with respect to the new probability
 measures depending on $n$.
Thus, while estimating \eqref{2.22.3} we assume that
$\pi_{t}\equiv0$.

Introduce
$$
x^{n}_{t} =
x^{\alpha^{n}_{\cdot}\hat{\bbeta} 
(\alpha^{n}_{\cdot}) x^{n}}_{t} ,\quad
 y^{n}_{t} 
=x^{p\alpha^{n}_{\cdot}\hat{\bbeta} 
( \alpha^{n}_{\cdot}) x^{n} }_{t} ,
$$
$$
c^{n}_{t}=c^{\alpha^{n}_{t}\hat{\beta}_{t}(\alpha^{n}_{\cdot})}
(x^{n}_{t}),\quad
pc^{n}_{t}=c^{p\alpha^{n}_{t}\hat{\beta}_{t}(\alpha^{n}_{\cdot})}
(y^{n}_{t})
$$
$$
f^{n}_{t}=f^{\alpha^{n}_{t}\hat{\beta}_{t}(\alpha^{n}_{\cdot})}
(x^{n}_{t}),\quad
pf^{n}_{t}=f^{p\alpha^{n}_{t}\hat{\beta}_{t}(\alpha^{n}_{\cdot})}
(y^{n}_{t})
$$
$$
r^{n}_{t}=r^{\alpha^{n}_{\cdot}\hat{\bbeta} 
( \alpha^{n}_{\cdot})}_{t},\quad
p\phi^{n}_{t}=\int_{0}^{t}[r^{n}_{s}]^{2}pc^{n}_{s}  \,ds,
$$
and define $\gamma^{n}$ as the first exit time of $y^{n}_{t}$
from $G$. Notice that, for any $T\in[0,\infty)$, \eqref{2.22.3} equals
$$
I_{1n}(T)+I_{2n}(T)-I_{3n}(T),
$$
where
$$
I_{1n}(T)=
E\int_{0}^{\tau^{n}\wedge\gamma^{n}\wedge T}
[r^{n}_{t}]^{2}\big[f^{n}_{t}\exp(-\phi^{n}_{t})
-pf^{n}_{t}\exp(-p\phi^{n}_{t})\big]\,dt,
$$
$$
I_{2n}(T)=E\int^{\tau^{n}}_{\tau^{n}\wedge\gamma^{n}\wedge T}[r^{n}_{t}]^{2}
f^{n}_{t}\exp(-\phi^{n}_{t})\,dt,
$$
$$
I_{3n}(T)=E\int^{\gamma^{n}}_{\tau^{n}\wedge\gamma^{n}\wedge T}[r^{n}_{t}]^{2}
pf^{n}_{t}\exp(-p\phi^{n}_{t})\,dt,
$$

By using the inequalities $|e^{-a}-e^{-b}|\leq|a-b|$
valid for $a,b\geq0$ and $|ab-cd|\leq |b|\cdot|a-c|+|c|\cdot
|b-d|$ and also using the boundedness of $r^{\alpha_{\cdot}\beta_{\cdot}}$,
$c^{\alpha\beta}$, and 
$f^{\alpha\beta}$ we easily conclude that
$$
|I_{1n}(T)|\leq N(1+T)E\int_{0}^{\tau^{n}\wedge T}\big[
|f^{n}_{t}-pf^{n}_{t}|+|c^{n}_{t}-pc^{n}_{t}|\big]\,dt .
$$
Observe that, if $\alpha^{n}_{t}\in A_{1}$, then
$$
|f^{n}_{t}-pf^{n}_{t}|\leq W_{f}(|x^{n}-y^{n}_{t}|),
$$
where $W_{f}$ is the modulus of continuity of $f^{\alpha\beta}(x)$
with respect to $x$
uniform with respect to $\alpha,\beta$. 
Similar estimate holds for $|c^{n}_{t}-pc^{n}_{t}|$ in which
$W_{c}$ is the modulus of continuity of $c^{\alpha\beta}(x)$.
Furthermore,
$$
E\int_{0}^{\tau^{n}\wedge T}I_{\alpha^{n}_{t}\in A_{2}}\,dt
\leq e^{T/\delta}E\int_{0}^{\tau^{n} }I_{\alpha^{n}_{t}\in A_{2}}
e^{-\phi^{n}_{t}}\,dt\leq Ne^{T/\delta}/K,
$$
where the last inequality is due to \eqref{2.22.4}.
Hence,
$$
|I_{1n}(T)|\leq N(1+T)^{2}E[W_{c}+W_{f}](
\sup_{t\leq \tau^{n}\wedge T} |x^{n}_{t}-y^{n}_{t}|)
+ Ne^{T/N}/K.
$$
We may and will assume that $W_{c}(r)$ and $W_{f}(r)$ are concave
functions on $[0,\infty)$, so that
$$
|I_{1n}(T)|\leq N(1+T)^{2}[W_{c}+W_{f}](E
\sup_{t\leq \tau^{n}\wedge T} |x^{n}_{t}-y^{n}_{t}|)
+ Ne^{T/N}/K.
$$
Next use the fact that as follows from Lemma \ref{lemma 3.24.1} 
$$
E\sup_{t\leq \tau^{n}\wedge T} |x^{n}_{t}-y^{n}_{t}|
\leq Ne^{NT}/\sqrt{K}.
$$
Then we conclude that
\begin{equation}
                                                               \label{2.22.5}
|I_{1n}(T)|\leq N(1+T)^{2}[W_{c}+W_{f}](Ne^{NT}/\sqrt{K})
+ Ne^{T/\delta}/K.
\end{equation}

While estimating $I_{2n}(T)$ we again use the boundedness
of the data and use Remark \ref{remark 2.22.2}
and by It\^o's formula to obtain that
$$
|I_{2n}(T)|\leq NEI_{\tau^{n}\geq\gamma^{n}\wedge T}
\int_{\gamma^{n}\wedge T}^{\tau^{n} }
[r^{n}_{t}]^{2}\,dt\leq NEI_{\tau^{n}\geq\gamma^{n}\wedge T}
G(x^{n}_{\gamma^{n}\wedge T})
$$
$$
\leq NEI_{\tau^{n}\geq\gamma^{n}\wedge T}|
G(x^{n}_{\gamma^{n}\wedge T})-G(y^{n}_{\gamma^{n}\wedge T})|
+NEI_{\tau^{n}\geq\gamma^{n}\wedge T}|G(y^{n}_{\gamma^{n}\wedge T})|
$$
$$
\leq NE\sup_{t\leq \tau^{n}\wedge T} |x^{n}_{t}-y^{n}_{t}|
+NEI_{\gamma^{n}>T}G(y^{n}_{T}).
$$
By Lemma 5.1 of \cite{Kr_13_1}
$$
EI_{\gamma^{n}>T}G(y^{n}_{T})\leq Ne^{-T/N}.
$$

Next,
$$
|I_{3n}(T)|\leq NEI_{\gamma^{n}\geq\tau^{n}\wedge T}
\int_{\tau^{n}\wedge T}^{\gamma^{n} }
[r^{n}_{t}]^{2}\,dt\leq NEI_{\gamma^{n}\geq\tau^{n}\wedge T}
G(y^{n}_{\tau^{n}\wedge T})
$$
$$
\leq
NE\sup_{t\leq \tau^{n}\wedge T} |x^{n}_{t}-y^{n}_{t}|
+NEI_{\gamma^{n}\geq\tau^{n}\wedge T}
G(x^{n}_{\tau^{n}\wedge T})
$$
$$
\leq NE\sup_{t\leq \tau^{n}\wedge T} |x^{n}_{t}-y^{n}_{t}|
+NEI_{\tau^{n}>\  T}
G(x^{n}_{ T}).
$$
We use again Lemma 5.1 of \cite{Kr_13_1} and conclude that, for $K\geq1$,
 \eqref{2.22.3} is less than
$$
 w(T,K):=N(1+T)^{2}[W_{c}+W_{f}](Ne^{N_{1}T}/\sqrt{K})
+  Ne^{N_{1}T}/\sqrt{K}+Ne^{-T/N_{2}}.
$$

Thus,
\eqref{2.22.6} yields
$$
v_{K}(x^{n})\leq 1/n
+v^{p\alpha_{\cdot}^{n}\bbeta(p\alpha_{\cdot}^{n})}(x^{n})
+w(T,K).
$$
Hence
$$
v_{K}(x^{n})\leq \sup_{\alpha_{\cdot}\in\frA}
 v^{ \alpha _{\cdot}
 \bbeta ( \alpha_{\cdot} )}(x^{n})+w(T,K) +1/n.
$$
Owing to  the arbitrariness of $\bbeta\in\bB$ we have
$$
v_{K}(x^{n})\leq v(x^{n})+w(T,K) +1/n,
$$
and  the arbitrariness of $x^{n}$ yields that for $K\geq1$
\begin{equation}
                                                           \label{2.23.2}
\sup_{\bar{D}}(v_{K}-v) \leq  w(T,K)  ,
\end{equation}
which leads to the desired result after first letting $K\to\infty$
and then $T\to\infty$. The theorem is proved.

\begin{remark}
                                                    \label{remark 2.23.1}
Assume that $c^{\alpha\beta}(x)$ and $f^{\alpha\beta}(x)$
are H\"older continuous with respect to $x$ with exponent $\kappa\in(0,1]$
and constant independent of $\alpha$ and $\beta$. Then by taking
$T$ such that $e^{N_{1}T}=K^{1/4}$ we see that, for $K\geq1$, the left-hand side of 
 \eqref{2.23.2} is dominated by
$$
N(1+\ln K)^{2}K^{-\kappa/4}+NK^{-1/(4N_{1}N_{2})}.
$$
Hence, there is a $\chi\in(0,1]$ such that 
the left-hand side of 
 \eqref{2.23.2} is dominated by $NK^{-\chi}$ for $K\geq 1$.
Thus, we have justified a claim made in Section 5 of \cite{Kr_14_1}.

\end{remark}

 \mysection{Proof of Theorem \protect\ref{theorem 1.14.1}}
                                                         \label{section 4.6.2}

The properties of $P$ listed before Theorem \ref{theorem 9.23.1}
or just the construction of $P$ in \cite{Kr_12}
yield that there is a set $A_{2}$, having no common points with $A$, and bounded 
continuous
functions $\sigma^{\alpha}=\sigma^{\alpha\beta}$,
$b^{\alpha}=b^{\alpha\beta}$, $c^{\alpha}=c^{\alpha\beta}$  
(independent of  $x$ and $\beta$), 
and $f^{\alpha\beta}\equiv0$ defined on $A_{2}$ 
such that  the assumptions in Section \ref{section 2.26.3}
are satisfied perhaps with   different constants $\delta $ and $K_{0}$  
and for $a^{\alpha}:=a^{\alpha\beta}
=(1/2)\sigma^{\alpha}(\sigma^{\alpha})^{*}$ we have
\begin{equation}
                                                  \label{2.8.5}
P[u](x)=\sup_{\alpha\in A_{2}}
\big[a_{ij}^{\alpha} D_{ij}u(x)
+b^{\alpha}_{i} D_{i}u(x)-c^{\alpha} u(x) \big].
\end{equation}

Use the notation from Section \ref{section 2.25.1} and observe that
$$
\max(H[u](x),
P[u](x)-K)
$$
$$
=\max\big\{\supinf_{\alpha\in A_{1}\,\,\beta\in B}
[L^{\alpha\beta}u(x)+f^{\alpha\beta}(x)],
\supinf_{\alpha\in A_{2}\,\,\beta\in B}
[L^{\alpha\beta}u(x)+f^{\alpha\beta}(x)-K]\big\}
$$
$$
=\supinf_{\alpha\in\hat{A}\,\,\beta\in B}
\big[ L^{\alpha\beta}u( x)+f^{\alpha\beta}_{K}( x)]
\quad (f^{\alpha\beta}_{K}( x)=f^{\alpha\beta} ( x)
I_{\alpha\in A_{1}}-KI_{\alpha\in A_{2}}),
$$
where the first equality follows from the definition of $H[u]$,
\eqref{2.8.5}, and the fact that $L^{\alpha \beta}$
is independent of $\beta$ for $\alpha\in A_{2}$.
It follows by Theorems \ref{theorem 9.23.1} and \ref{theorem 2.16.1}
that $u_{K}=v_{K}$ and by Theorem \ref{theorem 1.20.1} that in $\bar{G}$
$$
v=\lim_{K\to\infty}u_{K},
$$
where the right-hand side is indeed independent
of the probability space, filtration, and the choice of $w,r,\pi$.
Since the above convergence is uniform, $v$ is continuous   in $\bar{G}$.
The theorem is proved.

\end{document}